\documentclass[12pt,leqno]{amsart}
\usepackage[dvips]{graphicx}
\usepackage{amsfonts}
\usepackage{amssymb}
\usepackage{amsmath}
\usepackage{amscd}
\usepackage{graphicx}

\def\DATE{\today}
\newtheorem{theorem}{Theorem}
\newtheorem{definition}[theorem]{Definition}
\newtheorem{corollary}[theorem]{Corollary}

\newtheorem{lemma}[theorem]{Lemma}

\newtheorem{proposition}[theorem]{Proposition}

\newcommand\C{\mathbb{C}}
\newcommand\R{\mathbb{R}}
\newcommand\g{\mathfrak{g}}

\newcommand\h{\mathfrak{h}}
\newcommand\K{\mathbb{K}}
\newcommand\Z{\mathbb{Z}}
\newcommand\D{\delta_{\mu_0}}

\newcommand\HH{\mathbb{H}}

\newcommand\T{\mathcal{T}}

\newcommand\pf{\noindent{\it Proof. }}
\def\ra{\rightarrow}

\title{Cartan class of Invariant  forms on Lie groups.}
\author[Goze - Remm]{Michel Goze and Elisabeth Remm}
\address{M.G: Ramm Algebra Center, 4 rue de Cluny. F.68800 Rammersmatt}
\address{E.R: Laboratoire de Math\'ematiques et Applications,
        Universit\'e de Haute Alsace, Facult\'e des Sciences et
        Techniques, 4, rue des Fr\`eres Lumi\`ere,
        68093~Mulhouse~cedex, France.}
\email{Elisabeth.Remm@uha.fr, goze.RAC@gmail.com}

\begin{document}

\maketitle

\begin{abstract}
We are interested in the class, in the Elie Cartan sense, 
of left invariant forms on a Lie group. We construct the class of Lie algebras provided with a contact form and classify the frobeniusian Lie algebras up to contraction. We also study forms which are invariant by a subgroup. We show that the simple group $SL(2n,\R)$ which doesn't admit left invariant contact form, yet admits a contact form which is invariant by a maximal compact subgroup.

\end{abstract}

\section{Introduction}

The Elie Cartan class of a differential form $\omega$ on a differentiable manifold $M$ in a point $x$ of $ M$ corresponds 
to the number of Pfaffian forms linearly independent which are necessary to write the form $\omega$ and its differential in this point. It is a function on $M$ which is in general non constant. Here we consider manifolds $M$ which are Lie groups $G$  and the class of differential forms on $G$ which are either left invariant or invariant by a subgroup $H$ of $G$.  The left invariant forms are of constant class so we study two interesting  cases
- when $\dim G=2p+1$ and $\omega$ is of class $2p+1$ ($\omega$ is a contact form),
- when $\dim G=2p$ and $\omega$ is of class $2p$ ($\omega$ is exact and symplectic).
If $G$ is a Lie group provided with a left invariant contact form, its Lie algebra is a quadratic deformation of the Heisenberg algebra. We describe these deformations and, as an example, we focus on the case of nilpotent Lie algebras. If $G$ is a Lie group provided with an exact symplectic form, its Lie algebra is frobenusian. We give the classification up to contraction of these Lie algebras. This leads to a description of all the invariants up to isomorphism, invariant by deformation of these algebras. 

If the form $\omega$  on the group $G$ is invariant by a proper subgroup $H$ of $G$ the class of $\omega$ is not necessarily constant. If $G$ is provided with such a contact form the structure is described in \cite{Lutz}. After recalling  that a simple Lie group can not be provided with a left invariant contact form (or frobenusian) \cite{GozeCras1}, we describe a contact form on $SL(2n)$ invariant by the subgroup $SO(2n)$.

\tableofcontents

\section{Cartan class of a differential form}
\subsection{Definitions} 
Let $M$ be a $n$-dimensional differentiable manifold. We denote by $\T(M)$ the set of vector fields on $M$ and by $\Lambda^q(M)$ the set of  $q$-forms on $M$ (the vector fields and the exterior forms are always supposed to be differentiable). A $1$-form  is called a Pfaffian form. If $X \in \T(M)$ and $\theta \in \Lambda^p(M)$, the interior product of $X$ with $\theta$ is the $(p-1)$-form, $i(X)\theta$, defined by $i(X)\theta(X_1,\cdots,X_{p-1})=\theta(X,X_1,\cdots,X_{p-1})$ for any $X_1,\cdots,X_{p-1} \in \T(M)$.
The Lie derivative $L_X\theta$ of the $q$-form $\theta$ is the $p$-form $L_X\theta=d(i(X)\theta)+i(X)d\theta.$ 
\begin{definition}
Let $\theta$ be a $q$-form on the differential manifold $M$ with $q\leq n$. The Cartan class of $\theta$ at the point $x \in M$ is the codimension of the characteristic space $\mathcal{C}_x(\theta)$ defined by
$$\mathcal{C}_x(\theta)=\{X_x \in T_x(M), \ i(X_x)\theta(x)=0, \ i(X_x)(d \theta)(x)=0\}$$
where $T_xM$ is the tangent space to $M$ at the point $x$.
\end{definition}
The characteristic space $\mathcal{C}_x(\theta)$ is the intersection of  the associated subspaces 
$A_x(\theta)=\{X_x \in T_x(M), \ i(X_x)\theta(x)=0\}$ and $A_x(d\theta).$
In particular, if $\omega$ is a Pfaffian form, then
\begin{itemize}
\item $\omega$ is of Cartan class $2p+1$ at $x$ if $(\omega \wedge (d\omega)^p)(x) \neq 0$ and $   (d\omega)^{p+1}(x)=0$,
\item $\omega$ is of Cartan class $2p$ at $x$ if $(d\omega)^p(x) \neq 0 $ and $ \omega \wedge (d\omega)^{p}(x)=0$.
\end{itemize}
If $M=G$ is a connected Lie group, we denote by $\T_L(G)$ (respectively $\T_R(G)$) the set of left invariant vector fields on $G$ (resp. right invariant) and by  $\Lambda_L^p(G)$ the vector space of left invariant $p$-forms on $G$. The left invariance of $X \in \T(G)$ is equivalent to $L_{\widetilde{Y}}X=[X,\widetilde{Y}]=0$ for any $\widetilde{Y} \in \T_R(G)$. The vector space $\T_L(G)$ is identified with the Lie algebra $\g$ of $G$ and we will not distinguish the left invariant vector field $X$ on $G$ and the corresponding vector in $\g$, also denoted by $X$.
Similarly, a left invariant $p$-form on $G$ satisfies $L_{\widetilde{Y}}\theta=0$ for any $\widetilde{Y} \in \T_R(G)$. We identified the vector space $\Lambda_L^p(G)$ with the exterior space $\Lambda^p(\g)$ and more particulary any left invariant Pfaffian form $\omega \in \Lambda^1_L(G)$ with its  corresponding linear form, also denoted  by $\omega$, in the dual space $\g^*$. It is obvious that any left invariant form
$\omega \in \g^*$ is of constant class equal to the codimension of the characteristic subalgebra
$$\mathcal{C}(\omega)=\{X \in \g, \  \omega(X)=0, i(X)(d\omega)=0\}.$$
In this case $i(X)d(\omega)(Y)=-\omega[X,Y]$ for any $X,Y \in \g$ where $[X,Y]$ denotes the Lie bracket of $\g$.

\subsection{Contact and symplectic manifolds}

\begin{definition}
 If $M$ is $(2p+1)$-dimensional, a contact form is a Pfaffian form of Cartan class $2p+1$ at any point $x$ of $M$.

\noindent If $M$ is $(2p)$-dimensional, a symplectic form is a closed $2$-form $\theta$ of Cartan class $2p$ at any point $x$ of $M$.
\end{definition}
So we will speak about contact manifolds or symplectic manifolds and when $M=G$ is a Lie group, we will speak about contact or symplectic Lie group. If the contact form is left invariant, since it corresponds to an element of the dual of the Lie algebra $\g$ of $G$, we will say that $\g$ is a contact Lie algebra and the corresponding linear form on $\g$ a contact form on $\g$.

\subsection{Poisson brackets on contact manifolds}
Recall in a first step classical results on Poisson brackets defined on a symplectic manifold. Let $(M^{2p},\theta)$ be a symplectic manifold. For any Pfaffian form $\omega$ on $M$, there exists a unique vector field $X_\omega$ defined by $\omega=i(X_{\omega})\theta$. The Poisson bracket $\{\omega_1,\omega_2\}$ of the two Pfaffian forms $\omega_1,\omega_2$ which is associated with the symplectic form $\theta$ is the Pfaffian form $\{\omega_1,\omega_2\}=i([X_{\omega_1},X_{\omega_2}])\theta.$ Recall that if $\omega_1$ and $\omega_2$ are closed Pfaffian forms, then the Poisson bracket is an exact Pfaffian form equal to $\{\omega_1,\omega_2\}=-d(\theta(X_{\omega_1},X_{\omega_2})).$ Let $\mathcal{C}^\infty(M)$ be the associative algebra of differentiable functions on $M$. Then the bracket on $\mathcal{C}^\infty(M)$
$$\{f,g\}=-\theta(X_{df},X_{dg}), \  \forall f,g \in \mathcal{C}^\infty(M)$$
is a Lie bracket satisfying the Leibniz identity:
$$\{fg,h\}=f\{g,h\}+\{f,h\}g.$$
Thus $\mathcal{C}^\infty(M)$  is a Poisson algebra. This structure plays a fundamental role in the theory of deformations of associative algebras and in the theory of deformation  quantization.

\medskip

We can also define a Poisson bracket in the context of contact geometry.  Let $(M^{2p+1},\alpha)$ be a contact manifold, that is, $\alpha$ is a contact form on $M$. There exists one and only one vector field $Z_\alpha$ on $M^{2p+1}$, called the Reeb vector field of $\alpha$, satisfying the following properties:
\begin{itemize}
\item $\alpha(Z_\alpha)=1$,
\item $i(Z_\alpha)d\alpha=0$
\end{itemize}
at any point of $M^{2p+1}$. Let $\mathcal{D}_\alpha(M^{2p+1})$ be the set of first integrals of $Z_\alpha$, that is,
    $$\mathcal{D}_\alpha(M^{2p+1})=\{f \in \mathcal{D}(M^{2p+1}), \ Z_\alpha (f)=0\}.$$
    Since we have $Z_\alpha(f)=i(Z_\alpha)df=0$, thus $df$ is invariant by $Z_\alpha$.
    \begin{lemma}
    $\mathcal{D}_\alpha(M^{2p+1})$ is a  commutative associative subalgebra of $\mathcal{D}(M^{2p+1})$.
\end{lemma}
\pf In fact this is a consequence of the classical formulae:
$$Z_\alpha(f+g)=Z_\alpha(f)+Z_\alpha(g),  \ Z_\alpha(fg)=(Z_\alpha(f))g+f(Z_\alpha(g)).$$
\begin{lemma}
For any non zero Pfaffian form $\beta$ on $M^{2p+1}$ such that $\beta(Z_\alpha)=0$, there exists a vector field $X_\beta$ with $\beta(Y)=d\alpha (X_\beta,Y)$ for any vector field $Y$. Two vector fields $X_\beta$ and $X'_\beta$ satisfying this property are such that $i(X_\beta-X'_\beta)d\alpha=0$.
\end{lemma}
\pf For any point $x \in  M^{2p+1}$, there exists a Darboux open neighborhood $U$ of $x$ with a coordinates system $\{x_1,\cdots,x_{2p+1}\}$ such that the restriction  of the contact form to $U$  writes $\alpha=x_1dx_2+\cdots+x_{2p-1}dx_{2p}+dx_{2p+1}$. The Reeb vector field is $Z_\alpha=\partial/\partial x_{2p+1}$.  A Pfaffian form $\beta$ satisfying $\beta(Z_\alpha)=0$ writes $\beta=\sum_{i=1}^{2p}\beta_idx_i.$ We consider the vector field
$$\displaystyle X_\beta=\sum_{i=1}^{p}\left(\beta_{2i}\frac{\partial}{\partial x_{2i-1}}-\beta_{2i-1}\frac{\partial}{\partial x_{2i}}\right).$$
It satisfies
$$d\alpha (X_\beta,Y)=\sum\limits_{i=1}^p(dx_{2i-1}\wedge dx_{2i})(X_\beta,Y)=\sum\limits_{i=1}^p(\beta_{2i} dx_{2i}+\beta_{2i-1}dx_{2i-1})(Y)=\beta (Y).$$
If $X_\beta$ and $X'_\beta$ satisfy $\beta=i(X_\beta)d\alpha=i(X'_\beta)d\alpha$, then $i(X_\beta-X'_\beta)d\alpha=0$. The vector field $X_\beta-X'_\beta$ is in the distribution given by the associated subspaces $A(d\alpha(x))_{x \in M}.$  Thus, since $\alpha$ is a contact form,  for any vector field  $Y$ belonging to this distribution, we have $Y=aZ_\alpha$. We deduce that
$$X'_\beta=X_\beta + a Z_\beta.$$

Let  $f\in \mathcal{D}_\alpha(M^{2p+1})$, that is, satisfying $Z_\alpha (f)=0$.  Then the Pfaffian form $\beta=df$ satisfies $df(Z_\alpha)=0$. We can defined  a vector field  $X_{df}$ associated with $df$ and  we writes $X_f$ for $X_{df}$.
\begin{theorem}
The algebra  $\mathcal{D}_\alpha(M^{2p+1})$ is a Poisson algebra.
\end{theorem}
\pf Let $f_1,f_2$ be in $\mathcal{D}_\alpha(M^{2p+1})$. Since we have $d\alpha(X_{f_1},X_{f_2})=d\alpha(X_{f_1}+U_1,X_{f_2}+U_2)$ for any $U_1,U_2 \in A(d\alpha)$, the bracket
$$\{f_1,f_2\}_\alpha=d\alpha(X_{f_1},X_{f_2})$$
is well defined. In order to prove that it is a Poisson bracket, we have to verify that it is a Lie bracket satisfying the Leibniz identity. In the Darboux coordinates system, we have
$$\{f_1,f_2\}_\alpha=\sum_{i=1}^{p} \left(\frac{\partial f_2}{\partial x_{2i}}\frac{\partial f_1}{\partial x_{2i-1}}-\frac{\partial f_1}{\partial x_{2i}}\frac{\partial f_2}{\partial x_{2i-1}}\right).$$
Thus, in this system, we find  the classical expression of the canonical Poisson product in $\R^{2p}$ and Jacobi and Leibniz identities are satisfied.

\section{Invariant forms on Lie groups}
\subsection{$J$-invariant forms}
Let $G$ be a Lie group and $J$ a linear subspace of the Lie algebra $\widetilde{\g}=\T_R(G)$.
\begin{definition}A $p$-form $\theta$ on $G$ is $J$-invariant if $L_{\widetilde{Y}}\theta =0$ for any $\widetilde{Y} \in J$.
\end{definition}
In particular, if $ J=\widetilde{\g}$, a $J$-invariant form is left invariant.  In any point $g \in G$, $J$ defined a linear subspace $J_g$ of $T_gG$ and the distribution $\{J_g\}_{g \in G}$ is regular, that is, $\dim J_g=\dim J$. If $\omega$ is a $J$-invariant Pfaffian form on $G$, the singular set $\Sigma_\omega$ associated with $\omega$ is the subset of $G$:
$$\Sigma_\omega=\{ g \in G, \ \omega(\widetilde{Y})(g)=0 \ \forall \widetilde{Y} \in J\}.$$
\begin{lemma} (\cite{Lutz}). If $\omega$ is a $J$-invariant contact form on $G$, then for any basis $\{\widetilde{Y}_1,\cdots,\widetilde{Y}_k\} $ of $J$, the map
$$\varphi(g)=(\omega(\widetilde{Y}_1)(g),\cdots,\omega(\widetilde{Y}_k)(g))$$
with values in $\R^k$ is of rank $k$ on $\Sigma_\omega$ and of rank $k-1$ on $G - \Sigma_\omega$.
\end{lemma}
In particular, we deduce
\begin{itemize}
  \item If $k \leq p+1$, then $\Sigma_\omega=\emptyset.$
  \item If $J$ is an abelian subalgebra of $\widetilde{\g}$, then $k \leq p+1.$
\end{itemize}
For example, there exists on the torus $\mathbb{T}_5$ a $J$-invariant contact form where $J$ is the abelian $2$-dimensional Lie algebra \cite{Lutz}.

\subsection{A cohomology associated with $J$-invariant forms on a Lie algebra}

On a connected Lie group $G$, one defines classically the following complexes and their corresponding cohomologies:
\begin{enumerate}
  \item The de Rham cohomology of $G$, $H^*_{dR}(G)$, corresponding to the complex of differential forms on $G$.
  \item The cohomology of left invariant forms $H_L^*(G)$ corresponding to the complex of left invariant forms on $G$.
  \item The cohomology $H^*(G,M)$ of the group $G$ with values in a $G$-module $M$. The $p$-cochains are the maps on $G^p$ with values in $M$.
  \item The differentiable cohomology $H_{diff}^*(G,M)$ where $M$ is a differentiable $G$-module. The $p$-cochains are the differentiable maps on $G^p$ with values in $M$.
\end{enumerate}
The classical cohomologies associated with the Lie algebra $\g$ of $G$ are
\begin{enumerate}
  \item The Chevalley-Eilenberg cohomology $H_{CE}^*(\g,V)$ associated with the complex of exterior forms on $\g$ with values on a $\g$-module $V$.
  \item The $\h$-basic cohomology $H ^*(\g,\h,\K)$ where $\h$ is a Lie subalgebra of $\g$, and the $p$-cochains of the corresponding complex, the skew linear $p$-forms on $\g$ whose characteristic space contains  $\h$.
\end{enumerate}
  We have the following classical isomorphisms:
\begin{enumerate}
  \item $H_L^*(G)= H_{CE}^*(\g,\R)$
  \item If $G$ is compact, $H^*_{dR}(G)=H_L^*(G)$
  \item If $G$ is semi-simple and if $\h$ is the Lie algebra of the maximal compact subgroup of $G$, then $H^*_{diff}(G,\R)=H^*(\g,\h,\R)$.
\end{enumerate}
It is easy to define a complex on $G$ dependent of the notion of $J$-invariance.
Let $\Lambda^p_J(G)$ the space of $J$-invariant $p$-forms on the Lie group $G$. Since
$$L_X \circ d = d\circ L_X$$
for any vector field $X$, we have
$$\theta \in \Lambda^p_J(G) \Rightarrow d\theta \in \Lambda^{p+1}_J(G).$$
Then we define naturally a complex where the spaces of cochains are $\Lambda^p_J(G)$ and the coboundary operators are the classical exterior differential. We denote by $(\Lambda^p_J(G),d)_p$ this complex:
$$
\Lambda^0_J(G) \stackrel{d} \ra \Lambda^1_J(G) \stackrel{d} \ra \Lambda^2_J(G) \ldots \stackrel{d}\ra \Lambda^n _J(G) \stackrel{d}\ra 0$$
where $n=\dim G$ and $\Lambda^0_J(G)$ is the space of functions on $G$ which are $J$-invariant. The corresponding cohomology $H^*_J(G)$ will be called the $J$-invariant de Rham cohomology. If $J=\widetilde{\g}$, thus
$H^*_J(G)$ coincides with the Chevalley-Eilenberg cohomology  of $\g$ with values in the trivial module $\R$. If $G$ is compact, thus
$$H^*_{dR}(G)=H^*_{CE}(\g,\R)$$
where $ H^*_{dR}(G)$ is the de Rham cohomology of $G$. These equalities are in general not satisfied when $G$ is not compact. For example, if $G=\R^2$ and $J$ is the subalgebra generated by $\displaystyle \frac{\partial}{\partial x}$, then $\Lambda^0_J(G)=\{f:\R^2\rightarrow \R, \displaystyle \frac{\partial f}{\partial x}=0\}$, that is, $\Lambda^0_J(G)=\{f(y)\}$, but $\Lambda^0_{dR}(G)=\{f(x,y)\}$, $\Lambda^0_{\R^2}(G)=\R.$

\subsection{The Cartan class of left invariant forms}
We identify the set of left invariant Pfaffian forms on a Lie group with the dual $\g^ *$ of its Lie algebra $\g$. More generally, the set of left invariant exterior $p$-forms on $G$ is identified with $\Lambda^ p(\g^*).$  In this context, we can speak about the Cartan class of an element $\omega$  of $\g^*$ or  $\Lambda^ p(\g^*)$ and denote it by $cl(\omega)$.
\begin{proposition}
Assume that $\g$ is a finite dimensional nilpotent Lie algebra. Then the class of any non zero $\omega \in \g^*$ is always odd.
\end{proposition}
\pf This is a consequence of the Engel's theorem. Since $\g$ is nilpotent, for any $X \in \g$, $adX$ is a nilpotent operator. Let $\omega \neq 0$, $\omega \in \g^*$. If the Cartan class of $\omega$ is even, for example $cl(\omega)=2p$, we have $(d \omega)^p \neq 0$ and $\omega \wedge (d \omega)^p=0$. There exists a basis $\{X_1,\cdots,X_n\}$ of $\g$ such that its dual basis $\{\omega_1,\cdots,\omega_n\}$ satisfies
$$
\omega=\omega_1, \ \ d\omega_1=\omega_1 \wedge \omega_2 +\cdots+\omega_{2p-1} \wedge \omega_{2p}.$$
This implies that
$$[X_1,X_2]=-X_1+Y$$
with $Y \in \K\{X_2,\cdots,X_n\}$. Moreover, we have $[X_2,Y] \in \K\{X_2,\cdots,X_n\}$  from the canonical decomposition of $d\omega_1$. This implies that $1$ is an eigenvalue of $adX_2$ and $adX_2$ cannot be nilpotent. Then the Cartan class of any non trivial $\omega \in \g$ is not even.
\begin{corollary}\cite{GozeCras1}
A real or complex finite dimensional nilpotent Lie algebra is not a frobeniusian Lie algebra.
\end{corollary}
Recall that a frobeniusian Lie algebra is of even dimension $n=2p$ and admits, by definition, a linear form  $\omega$ of class $2p$. Since the Cartan class of a nilpotent Lie algebra is always odd, we deduce the Corollary.

\medskip

 Let us note that this property is not true for solvable not nilpotent Lie algebras. But we have also:
 \begin{proposition}
 Let $\g$ be a real compact Lie algebra. For any $\omega \in \g^*$, $\omega \neq 0$, then its Cartan class is odd.
 \end{proposition}
 \pf Recall that a real Lie algebra $\g$ is called compact if there exists a compact Lie group with Lie algebra isomorphic to $\g$. We consider the compact Lie group $G$ whose Lie algebra is $\g$. Assume that $\omega$ is of even class, that is, there exists $k$ with $(d\omega)^k \neq 0$ and $\omega  \wedge (d\omega)^k =0.$ Then the exterior $(2k-1)$-form $\omega  \wedge (d\omega)^{k-1}$ is not closed because $d(\omega  \wedge (d\omega)^{k-1})=(d\omega)^k \neq 0$. This is impossible from the Stokes theorem on $G$ because $G$ is compact.

Since a compact Lie algebra is reductive, that is, a direct sum of a semi-simple Lie algebra and an abelian ideal, this leads to  study  the Cartan class of semi-simple Lie algebras.
\begin{proposition}\cite{BW,GozeCras1}\label{simple}
Let $\g$ be a real or complex semi-simple Lie algebra of rank $r$. Then, for any $\omega \in \g^*$ we have
$$cl(\omega) \leq n-r+1.$$
\end{proposition}
\pf For any $\omega  \in g^*$ with $\omega  \neq 0$, the associated space $A(d\omega)=\{X \in \g, i(X)d\omega=0\}$ is a subalgebra of dimension less than or equal to $n-c+1$ where $n=\dim \g$ and $c$ is the Cartan class of $\omega$. If $K$ denotes the Killing-Cartan form on $\g$, since $K$ is non degenerate, there exists $U_\omega$ such that $K(U_\omega,Y)=\omega(Y)$ for any $Y \in \g$. Since $K$ is $ad(\g)$-invariant, $A(d\omega)$ coincides with the commutator of $U_\omega$, that is, $\{X \in \g, \ [X,U_\omega]=0\}.$ But the dimension of this subspace is greater than or equal to $r$. Then $n-c+1 \geq r$, that is, $c \leq n-r+1.$
\begin{corollary}
Let $\g$ be a semi-simple Lie algebra of rank $r$. There exists a form of maximal class $c=\dim \g$ if and only if $r=1$. In particular none of semi-simple Lie algebra is frobeniusian.
\end{corollary}
\pf From the previous inequality, $c=n$ implies $r=1$. Thus $\g$ is of rank $1$ and it is isomorphic to $sl(2)$ in the complex case,  $sl(2,\R)$ or $so(3)$ in the real case. It is easy in each case to find a contact form. These algebras are of dimension $3$, so we deduce that $\g$ is not frobeniusian.
\begin{proposition}
Let $\g$ be a real $n$-dimensional Lie algebra such that for any $\omega  \in g^*$ with $\omega \neq 0$ we have $cl(\omega)=n$. Then $\g$ is of dimension $3$ and is isomorphic to $so(3)$.
\end{proposition}
\pf Let $I$ be a non trivial abelian ideal of $\g$. There exists $\omega \neq 0 \in \g^*$ such that $\omega (X)=0$ for any $X \in I$. This implies that $I \subset A(d\omega)$ and $cl(\omega) \leq n - dim I$ and this form is not of class $n$. Thus $\g$ is semi-simple. But a semi-simple Lie algebra with a form of maximal class is of rank $1$. It is isomorphic to $sl(2)$ or $so(3)$. In the first case, we have a form of class $2$. In the second case, we have a basis of $so(3)^*$ such that $d\omega_1=\omega_2\wedge \omega_3$, $d\omega_2=\omega_3\wedge \omega_1$, $d\omega_3=\omega_1\wedge \omega_2$. In this case any $\omega\neq 0 \in \g^*$ is of class $3$.
\medskip

\noindent{\bf Remarks}
\begin{enumerate}
\item {\bf On the class of linear form of simple Lie algebras.} In Proposition \ref{simple} we give an upper bound of the set of classes of linear forms on a semi-simple Lie algebra. It is easy to give a non null form which realizes this bound. But, in \cite{Goze3cycle}, we gave also a lower bound when  $\g$ is complex, simple and classical:
{\it The Cartan class of any linear non trivial form on a simple non exceptional Lie algebra of rank $r$ satisfies
$$cl(\omega) \geq 2r.$$}
\item{\bf Cartan class and the index of a Lie algebra.} For any $\omega \in \g^*$, we consider the stabilizer $\g_{\omega }=\{X\in \g, \omega\circ ad X=0\}$ and we consider $d$ the minimal dimension of $\g_{\omega }$ when $\omega$ lies in $\g^*.$ If $\omega$ is such that $\dim \g_{\omega }=d$, then, from \cite{Vergne1}, $\g_{\omega }$ is an abelian subalgebra of $\g$. From the point of view of the Cartan class, $\g_{\omega }$ is the characteristic space associated with $d\omega$ and the minimality is realized by a form of maximal class and we have $d=n-cl(\omega)+1$ if the Cartan class $cl(\omega)$ is odd or $d=n-cl(\omega)$ if  $cl(\omega)$ is even. In particular, if $\g$ is frobeniusian, then the maximal class is $2p=n$ and $d=0$. If $\g$ is provided with a contact form, then $d=n-n+1=1.$ For example, if $\g$ is a naturally graded filiform Lie algebra, it is isomorphic to $L_n$ or $Q_{2p}$ \cite{GozeKhakimbook}. In the first case the $cl(\omega) \in \{1,3\}$ and $d=n-2$, in the second case $cl(\omega) \in \{1,3,2p-1\}$ and $d=2$.
\end{enumerate}

\section{Lie groups with a left invariant contact form}

\subsection{A model up to contraction of contact Lie algebra}

Let $G$ be a $(2p+1)$-dimensional connected Lie group with Lie algebra $\g$. Recall that $\Lambda^1_L(G)=\Lambda^1_{J=\widetilde{\g}}(G)=\g^*.$ If $\omega \in \Lambda^1_L(G)=\g^*$ is a contact form, that is, $\omega \wedge (d\omega) ^p \neq 0$, we will say that the Lie algebra $\g$ admits a contact form. The simplest example of Lie algebra with a contact form is the Heisenberg algebra $\h_{2p+1}$. There exists a basis $\{X_1,\cdots, X_{2p+1}\}$ of $\h_{2p+1}$ such that $[X_{2k-1},X_ {2k}]= X_{2p+1}$ for $k=1,\cdots,p.$ If $\{\omega_1,\cdots, \omega_{2p+1}\}$ is the dual basis, the Maurer Cartan equations related to this basis are $d\omega_i=0, \ i=1, \cdots, 2p, \ d\omega_{2p+1}=-\omega_1\wedge\omega_2 -\cdots -\omega_{2p-1}\wedge\omega_{2p}.$ Thus $\omega_{2p+1}$ is a contact form on $\h_{2p+1}$.
\begin{theorem}
Any $(2p+1)$-dimensional Lie algebra $\g$ with a contact form is isomorphic to a quadratic formal deformation of $\h_{2p+1}$.
\end{theorem}
\pf A formal deformation $\g_t$ of a Lie algebra $\g_0$ is a Lie algebra whose  Lie bracket $\mu_t$ satisfies
$$\mu_t(X,Y)=\mu_0(X,Y)+t\varphi_1(X,Y)+t^2\varphi_2(X,Y)+\cdots +t^n\varphi_n(X,Y)+ \cdots $$
where $\mu_0(X,Y)$ is the Lie bracket of $\g_0$ and the maps $\varphi_i$ are bilinear on $\g_0$ with values in $\g_0$ (the Lie algebras $\g_0$ and $\g_t$ have the same underlying vector space). The formal deformation is called quadratic if $\varphi_i=0$ for any $i \geq 3$. It is called linear if $\varphi_i=0$ for $i \geq 2.$ If $\g_t$ is a quadratic formal deformation of $\g_0$, then the Jacobi conditions related to the bracket of $\g_t$ are equivalent to
\begin{equation} \label{quadra}
\left\{
\begin{array}{l}
  \delta_{\mu_0} \varphi_1 =0, \\
 \varphi_1\circ \varphi _1+ \delta_{\mu_0}  \varphi_2=0, \\
  \varphi_1\circ \varphi _2+\varphi_2\circ \varphi _1=0, \\
  \varphi_2\circ \varphi _2=0,
\end{array}
\right.
\end{equation}
where $\delta_{\mu_0} $  is the coboundary operator of the Chevalley-Eilenberg cohomology of $\g$ with values in $\g$, and if $\varphi$ and $\psi$ are bilinear maps, then $\varphi \circ\psi$ is the trilinear map given by
$$\varphi \circ\psi(X,Y,Z)=\varphi(\psi(X,Y),Z)+ \varphi(\psi(Y,Z),X)+\varphi(\psi(Z,X),Y).$$ In particular $\varphi\circ\varphi=0$ is equivalent to the Jacobi identity and $\varphi$, in this case, is a Lie bracket.

Let $\g$ be a $(2p+1)$-dimensional Lie algebra with a contact form $\omega$. There exists  a basis $\{\omega_1,\cdots, \omega_{2p+1}\}$ of $\g^*$ such that
$$
\left\{
\begin{array}{l}
\omega=\omega_{2p+1},\\
d\omega_{2p+1}=-\omega_1\wedge\omega_2 -\cdots -\omega_{2p}\wedge\omega_{2p+1}.
\end{array}
\right.
$$
We denote by $\{X_1,\cdots, X_{2p+1}\}$ the dual basis of $\g$. The constant structures of $\g$ in this basis are:

\begin{equation}\label{basis}
\left\{
\begin{array}{l}
\medskip
\mu(X_{2k-1},X_{2k})=X_{2p+1}+ \sum_{s=1}^{2p}C_{2k-1,2k}^sX_s, \ k=1,\cdots,p,\\
\medskip
\mu( X_l,X_r )=\sum_{s=1}^{2p}C_{l,r}^sX_s, \ 1\leq l<r\leq 2p+1, (l,r)\neq (2k-1,2k).
\end{array}
\right.
\end{equation}
 Let $f_t$  be the linear isomorphism of $\g$ given by $f_t(X_i)=tX_i$ for $i=1,\cdots,2p$ and $f_t(X_{2p+1})=t^2X_{2p+1}$ with $t \neq 0.$ Let $\g_t$ be the Lie algebra isomorphic to $\g$ defined by $f_t$. Its Lie bracket satisfies
$$\left\{
\begin{array}{l}
\medskip
\mu_t(X_{2k-1},X_{2k})=X_{2p+1}+ t\sum_{s=1}^{2p}C_{2k-1,2k}^sX_s, \ k=1,\cdots,p,\\
\medskip
\mu_t(X_l,X_r)=t\sum_{s=1}^{2p}C_{l,r}^sX_s, \ 1\leq l<r\leq 2p, (l,r)\neq (2k-1,2k),\\
\medskip
\mu_t(X_l,X_{2p+1})=t^2\sum_{s=1}^{2p}C_{l,{2p+1}}^sX_s, \ l=1,\cdots,2p.
\end{array}
\right.
$$
Let $\varphi_1$ and $\varphi_2$ be the skew symmetric bilinear maps
\begin{equation}\label{phi1}
\left\{
\begin{array}{l}
\medskip
\varphi_1(X_{2k-1},X_{2k})=\sum_{s=1}^{2p}C_{2k-1,2k}^sX_s, \ k=1,\cdots,p,\\
\medskip
\varphi_1(X_l,X_r)=\sum_{s=1}^{2p}C_{l,r}^sX_s, \ 1\leq l<r\leq 2p, (l,r)\neq (2k-1,2k),\\
\end{array}
\right.
\end{equation}
and
\begin{equation}\label{phi2}
\left\{
\begin{array}{l}
\medskip
\varphi_2(X_l,X_r)=0, \ \ l,r=1,\cdots,2p,\\
\medskip
\varphi_2(X_l,X_{2p+1})=\sum_{s=1}^{2p}C_{l,{2p+1}}^sX_s \ l=1,\cdots,2p,
\end{array}
\right.
\end{equation}
and the non defined values are equal to $0$. Then the Lie bracket $\mu$ of $\g$ is isomorphic to the quadratic deformation
$$\mu_t(X,Y)=\mu_0(X,Y)+t\varphi_1(X,Y)+t^2\varphi_2(X,Y)$$ where $\mu_0$ is the Lie bracket of the Heisenberg algebra $\h_{2p+1}$. $\clubsuit$

\medskip

\subsection{ Lie algebras with a contact form defined by a linear deformation}

 Let $\g$ be such a Lie algebra.  Its Lie bracket can be written:
$$\mu=\mu_0 +t\varphi_1$$
where $\mu_0$ is the Lie multiplication of the Heisenberg algebra and the skew symmetric bilinear map $\varphi_1$ satisfies the relations
$$
\left\{
\begin{array}{l}
\delta_{\mu_0}\varphi_1=0,\\
\varphi_1\circ \varphi_1=0,
\end{array}
\right.
$$
and Relations (\ref{phi1}) in the classical basis of the Heisenberg algebra. Here, the operator $\delta_{\mu_0}$ is the coboundary operator associated with the Chevalley-Eilenberg cohomology of the Heisenberg algebra. 

If $\varphi_1=\delta_{\mu_0}  f$ is a coboundary then the linear deformation is trivial, that is, isomorphic to the Heisenberg algebra (see \cite {Gerstenhaber}, \cite{Fialowski}, \cite{GozeNato}).
Assume that $\varphi_1$ is not a coboundary. We have
$$\D \varphi_1(X_{2k-1},X_{2k},X_{i})=\varphi_1(X_{2p+1},X_{i})+\lambda_i X_{2p+1}.$$
Since $\D \varphi_1=0$, thus  $\varphi_1(X_{2p+1},X_{i})=-\lambda_i X_{2p+1}.$ Thus we can find a linear endomorphism $f$ such that $(\varphi_1-\D f)(X_{2p+1},X_{i})=0.$ From the first step, the cocycle $\varphi_1$ can be chosen up to a coboundary. Then we can assume that $\varphi_1$ satisfies Relations (\ref{phi1}) and
$$\varphi_1(X_{i},X_{2p+1})=0, \ i=1,\cdots, 2p.$$
\begin{theorem}\label{symplec}
Any $(2p+1)$-dimensional Lie algebra $\g$ isomorphic to a linear deformation of the Heisenberg algebra $\h_{2p+1}$ is a central extension of a $2p$-dimensional symplectic Lie algebra  by its symplectic form.
\end{theorem}
\pf Let $\frak{k}$ be the $2p$-dimensional vector space generated by $\{X_1,\cdots,X_{2p}\}$. The restriction to $\frak{k}$ of the $2$-cocycle $\varphi_1$ is in values in $\frak{k}$. Since $\varphi_1 \circ \varphi_1=0$, it defines on
$\frak{k}$ a structure of $2p$-dimensional Lie algebra. If $\{\omega_1,\cdots,\omega_{2p+1}\}$ is the dual basis of the given classical basis of $\h_{2p+1}$, thus the condition $\D\varphi_1=0$ is equivalent to $d(\omega_1\wedge\omega_2+\cdots +\omega_{2p-1}\wedge\omega_{2p})=0.$ It implies that $\theta =\omega_1\wedge\omega_2+\cdots +\omega_{2p-1}\wedge\omega_{2p}$ is a closed $2$-form on $\frak{k}$ and $\g$ is a central extension of $\frak{k}$ by $\theta$. $\clubsuit$

\subsection{ Lie algebras with a contact form defined by a quadratic deformation}
We want to determine  $(2p+1)$-dimensional Lie algebras whose Lie bracket $\mu$ is a quadratic deformation of $\mu_0$:
$$\mu(X,Y)=\mu_0(X,Y) +t\varphi_1(X,Y)+t^2\varphi_2(X,Y)$$
where the bilinear form $\varphi_1$ and $\varphi_2$ satisfies Relations  (\ref{phi1}),  (\ref{phi2}) and  (\ref{quadra}). Remark that Relations (\ref{phi2}) implies $\varphi_2\circ \varphi _2=0.$ Thus $ \varphi _2 $ is a
Lie bracket. We denote by $\g_2$ the Lie algebra whose Lie bracket is $\varphi _2$. To simplify put $a_l^s=C_{l,2p+1}^s$ and Relation  (\ref{phi2}) becomes
$$\varphi_2(X_l,X_{2p+1})=\sum_{s=1}^{2p}a_l^sX_s.$$ Let $\frak{k}$ be the $2p$-dimensional vector space generated by $\{X_1,\cdots,X_{2p}\}$ and $f$  
the endomorphism  of $\frak{k}$ defined by $f(X_l)=\sum_{i=1}^{2p}a_l^sX_s$ for any $l=1,\cdots,2p.$ 
Then the Lie algebra $\g_2$ is an extension of the $2p$-dimensional abelian Lie algebra by $f$.

Consider the equation $ \varphi_1\circ \varphi _1+ \delta_{\mu_0}  \varphi_2=0.$
We have 
$$
\left\{
\begin{array}{l}
\medskip
\D \varphi_2(X_{2k-1},X_{2k},X_{l})=-f(X_l), \ i=l,\cdots,2p,  l \neq 2k-1,2k,\\
\medskip
\D \varphi_2(X_{l},X_{s},X_{2p+1})=\mu_0(f(X_s),X_l)-\mu_0(f(X_l),X_s).
\end{array}
\right.
$$
Since the skew symmetric map $\varphi_1$ satisfies Relation (\ref{phi1})
$$\varphi_1\circ \varphi _1(X_{l},X_{s},X_{2p+1})=0.$$
Then 
\begin{equation}\label{f}
\mu_0(f(X_i),X_j)=-\mu_0(X_i,f(X_j)), \ 1 \leq i < j \leq 2p.
\end{equation}
\begin{lemma}
The vector space
$$F=\{f \in End(\frak{t}),  \mu_0(f(X_i),X_j)=-\mu_0(X_i,f(X_j))\}$$
is a Lie subalgebra of $sl(2p+1, \R)$ isomorphic to $so(2p+1)$.
\end{lemma}
\pf If $f,g \in F$, then
$$\mu_0(g\circ f(X_i),X_j)=-\mu_0(f(X_i),g(X_j))=\mu_0(X_i,f\circ g(X_j)).$$
Then
$$\mu_0((g\circ f-f\circ g)(X_i),X_j)+\mu_0(X_i,(g\circ f-f\circ g)(X_j))=0$$
and $g\circ f-f\circ g \in F$, this shows that $F$ is a Lie subalgebra of $End(\frak{t})$. The matrix $(a_i^j)$ satisfies the relations 
$$
\left\{
\begin{array}{l}
a_{2k-1}^{2s}-a_{2s-1}^{2k}=0,\\
a_{2k-1}^{2s-1}+a_{2s}^{2k}=0,\\
a_{2k}^{2s}+a_{2s-1}^{2k-1}=0,\\
a_{2k}^{2s-1}-a_{2s}^{2k-1}=0.\\
\end{array}
\right.
$$
This implies that $\dim F=p(2p+1)$ and for $s=k$, $a_{2k-1,2k-1}+a_{2k,2k}=0.$ Thus the trace of $f$ is $0$ and $F$ is a Lie subalgebra of $sl(2p+1,\R)$.

\medskip

\noindent{\bf Remarks. } 
\begin{enumerate}
\item If $\delta_{\mu_0}\varphi_2=0$, then  $\varphi_1$ is a Lie bracket. Moreover, if $p > 1$, that is,  $\dim\g > 3$, 
$$f(X_i)=-\delta_{\mu_0}\varphi_2(X_{2k-1},X_{2k},X_i)=0$$
with $\ k=1,\cdots,p$  and $i=1,\cdots,2p, \ i \neq 2k-1$ or $2k$. Then $\delta_{\mu_0}\varphi_2 =0$ is equivalent to $f=0$ and it is equivalent to $\varphi_2=0$. The quadratic deformation is linear. The case $\dim \g=3$ will be studied in the next section.
\item If  $\delta_{\mu_0}\varphi_2 \neq 0$, the bilinear map $\varphi_1$ is a $2$-cocycle for the Heisenberg algebra and satisfies
$$\varphi_1\circ \varphi _1=-\D \varphi_2.$$
Then $\varphi_1$ is not a Lie bracket because Jacobi identity is satisfies up to a $2$-cocycle. This means that $\varphi_1$ defines a structure of $L_3$ non graded algebra (or a $L_{\infty}$ non graded algebra, see \cite{MarklLada} for the definitions). Although $\varphi_1$ does not define a Lie  algebra structure, we can define a notion of differential. If $\varphi_1(X_i,X_j)=\sum A_{ij}^kX_k$, we put
$$\displaystyle d_{\varphi_1}\omega_i=\sum_{1 \leq j< k \leq 2p+1} A_{jk}^i\omega _j \wedge \omega_k$$
and $$d_{\varphi_1}(\omega_i \wedge \omega_j)=d_{\varphi_1}(\omega_i )\wedge \omega_j-\omega_i \wedge d_{\varphi_1}(\omega_j).$$
If $\varphi_1$ satisfies Relation (\ref{phi1}), it induces a skew bilinear map $\tilde{\varphi_1}$ on $\frak{k}$, the linear space generated by $\{X_1,\cdots,X_{2p}\}$,  and the (non Lie) algebra $(\frak{k},\tilde{\varphi_1})$ is a differential  $L_{\infty}$ non graded algebra of dimension $2p$. 
\end{enumerate}
From these remarks, we deduce
\begin{proposition}
The skew bilinear map $\varphi_1$ satisfying Relation (\ref{phi1}) is a $2$-cocycle for the Heisenberg algebra $\h_{2p+1}$ if and only if the exterior $2$-form $\theta=\omega_1\wedge \omega_2+\cdots+\omega_{2p-1}\wedge \omega_{2p}$ is symplectic in $(\frak{k},\tilde{\varphi_1})$, that is, $d_{\tilde{\varphi_1}}\theta=0.$
\end{proposition}

\medskip

Consider the equation  
\begin{equation}\label{phi1phi2}
\varphi_1\circ\varphi_2+\varphi_2\circ\varphi_1=0.
\end{equation}
We have $(\varphi_1\circ\varphi_2+\varphi_2\circ\varphi_1)(X_i,X_j,X_k)=0$ when $1 \leq i,j,k,\leq 2p.$ Then (\ref{phi1phi2}) is satisfied if and only if 
$$(\varphi_1\circ\varphi_2+\varphi_2\circ\varphi_1)(X_i,X_j,X_{2p+1})=0$$
for any $1\leq i < j \leq 2p.$  But
$$
\begin{array}{lll}
(\varphi_1\circ\varphi_2+\varphi_2\circ\varphi_1)(X_i,X_j,X_{2p+1}) & =& \varphi_1(\varphi_2(X_j,X_{2p+1}),X_i)+\varphi_1(\varphi_2(X_{2p+1},X_i),X_j)\\
&&+\varphi_2(\varphi_1(X_i,X_j),X_{2p+1})\\
&= & \varphi_1(f(X_j),X_i)+ \varphi_1(X_j,f(X_i))-f( \varphi_1(X_j,X_i)),\\
\end{array}
$$
that is, $f$ is a derivation of the $2p$-dimensional algebra whose multiplication is $\varphi_1$.
Since 
$$\varphi_2(f(X_i),X_{2p+1})+ \varphi_2(X_i,f(X_{2p+1}))-f( \varphi_2(X_i,X_{2p+1}))=f(f(X_i))-f(f(X_i))=0,$$
$f$ is also a derivation of the Lie algebra $\g_2$. 
\begin{theorem}
Let $\g$ be a $(2p+1)$-dimensional Lie algebra.  Then  its bracket $\mu$ is isomorphic to a quadratic deformation $\mu_0+t\varphi_1+t^2\varphi_2$ of the Heisenberg algebra $\h_{2p+1}$ if and only if
\begin{enumerate}
\item There exists a graded derivation $f$ of $\h_{2p+1}$ associated to the $\Z_2$-grading $\h_{2p+1}=\frak{k} \oplus Z(\h_{2p+1})$ which is trivial on the center $Z(\h_{2p+1})$ such that $\varphi_2$ is a Lie bracket of a $(2p+1)$-dimensional Lie algebra obtained by an extension by $f$ of the abelian Lie algebra $\frak{k}$,
\item $\varphi_1$ is a multiplication of a non graded differential $L_3$-algebra on $\frak{k}$ associated to the triple $(\mu_0,0,\varphi_2)$ of bilinear map, admitting a symplectic form.
\end{enumerate}
\end{theorem}

\medskip

We are going to describe the quadratic deformations of the Heisenberg algebra  $\h_{2p+1}$ when $f$ is a diagonal derivation and $p \geq 2$. Assume that we have $f(X_i)=\rho_i X_i$ and $f(X_{2p+1})=0.$ Since $f$ is a derivation of $\h_{2p+1}$, the eigenvalues satisfy
$$\rho_1+\rho_2=\rho_3+\rho_4=\cdots=\rho_{2p-1}+\rho_{2p}=0.$$
We put $\lambda_k=\rho_{2k-1}$, that is,
$$f(X_{2k-1})=\lambda_k X_{2k-1}, \ f(X_{2k})=-\lambda_k X_{2k}, \ \ k=1,\cdots,p.$$
Thus we can assume that all the $\lambda_i$ are nonnegative and more precisely
$$\lambda_1 \geq \lambda_2 \geq \cdots \geq \lambda_p \geq 0.$$
If $A_{ij}^k$ are the structure constants of $\varphi_1$, the identity $\varphi_1\circ\varphi_2+\varphi_2\circ\varphi_1=0$ is equivalent to
$$\begin{array}{ll}
    \left\{
    \begin{array}{l}
    (\lambda_s-\lambda_k-\lambda_l)A_{2k-1,2l-1}^{2s-1}=0,\\
    (-\lambda_s-\lambda_k-\lambda_l)A_{2k-1,2l-1}^{2s}=0,
    \end{array}
    \right. & \left\{
    \begin{array}{l}
    (\lambda_s-\lambda_k+\lambda_l)A_{2k-1,2l}^{2s-1}=0,\\
    (-\lambda_s-\lambda_k+\lambda_l)A_{2k-1,2l}^{2s}=0,
    \end{array}
    \right. \\
    \left\{
    \begin{array}{l}
    (\lambda_s+\lambda_k-\lambda_l)A_{2k,2l-1}^{2s-1}=0,\\
    (-\lambda_s+\lambda_k-\lambda_l)A_{2k,2l-1}^{2s}=0,
    \end{array}
    \right. & \left\{
    \begin{array}{l}
    (\lambda_s+\lambda_k+\lambda_l)A_{2k,2l}^{2s-1}=0,\\
    (-\lambda_s+\lambda_k+\lambda_l)A_{2k,2l}^{2s}=0.
    \end{array}
    \right.
  \end{array}
$$
Since
$$
\left\{
\begin{array}{l}
\delta_{\mu_0}\varphi_2(X_{2k-1},X_{2k},X_l)=-\rho_lX_l, \ k=1,\cdots,p,\\
\delta_{\mu_0}\varphi_2(X_{i},X_{j},X_l)=0, \   i \leq j \leq k, (i,j)  \neq (2k-1,2k), (j,l) \neq (2k-1,2k). \\
\end{array}
\right.
$$
Then 
$$
\left\{
\begin{array}{l}
\varphi_1 \circ \varphi_1(X_{2k-1},X_{2k},X_l)=-\rho_lX_l, \ k=1,\cdots,p,\\
\varphi_1 \circ \varphi_1(X_{i},X_{j},X_l)=0, \   i \leq j \leq k, (i,j)  \neq (2k-1,2k), (j,l) \neq (2k-1,2k) .\\
\end{array}
\right.
$$
Since $f$ is also a derivation of $\varphi_1$, we have
$$(\rho_i+\rho_j)\varphi_1(X_i,X_j)=f(\varphi_1(X_i,X_j))$$
and $\varphi_1(X_i,X_j)=0$ as soon as $\rho_i+\rho_j$ is not a root of $f$. In particular 
$$\varphi_1(X_1,X_{2k-1})=\varphi_1(X_2,X_{2k})=0$$
for any $k$. Assume that $f$ is non singular, that is, $\lambda_i \neq 0$ for $i=1,\cdots,p.$   This implies
$f(\varphi_1(X_{2k-1},X_{2k}))=0$ and
$$ \varphi_1(X_{2k-1},X_{2k})=0, \ 1 \leq k \leq p.$$
We have
$$\varphi_1 \circ \varphi_1(X_{1},X_{2},X_{2k})=-\lambda_kX_{2k}=\varphi_1(\varphi_1(X_{2k},X_1),X_2)$$
and $\varphi_1(X_1,X_{2k})\neq 0.$
Thus $f(\varphi_1(X_1,X_{2k}))=(\lambda_1-\lambda_k)\varphi_1(X_1,X_{2k}) \neq 0$
and $\lambda_1-\lambda_k$ is a root for any $k >1.$ We deduce that the positive roots satisfy
$$\lambda_1- \lambda_k=\lambda_{p-k+2}, \  k = 2, \cdots, p.$$
Similarly, for any $k$ such that $2 < k <p$, 
$$\varphi_1 \circ \varphi_1(X_{3},X_{4},X_{2k})=-\lambda_kX_{2k}=\varphi_1(\varphi_1(X_{2k},X_3),X_4)$$
and $\lambda_2-\lambda_k$ is a positive root. If $\lambda_2-\lambda_{p-1}=\lambda_3$ then $\lambda_1=\lambda_2$ and $\lambda_p=0.$ This contradicts the hypothesis on $f$. Then  $\lambda_2-\lambda_{p-1}=\lambda_4$. This gives
$$\lambda_2-\lambda_{p-1-k}=\lambda_{k+4}, \ k=0,\cdots, p-4.$$
We deduce
$$\lambda_1-\lambda_2=\lambda_2-\lambda_3=\lambda_3-\lambda_4=\cdots=\lambda_{p-2}-\lambda_{p-1}=\lambda_p.$$
In particular
$$\lambda_k=\lambda_1-(k-1)\lambda_1$$
and
$$\lambda_{2}-\lambda_{p-1}=(p-3)\lambda_1=\lambda_4=\lambda_1-3\lambda_p.$$
This implies
$$\lambda_1=p\lambda_p; \ \lambda_2=(p-1)\lambda_p; \ \cdots ; \lambda_{p-1}=2\lambda_p,$$ 
that is,
$$f(X_{2k-1})=(p-k+1)\lambda_pX_{2k+1}, \ f(X_{2k})=-(p-k+1)\lambda_p X_{2k}.$$
We deduce
$$\varphi_2(X_{2k-1},X_{2p+1})=(p-k+1)\lambda_pX_{2k-1}, \ \varphi_2(X_{2k},X_{2p+1})=-(p-k+1)\lambda_pX_{2k},$$
and
$\varphi_1$ is $\Z$-graded by the roots of $f$.

\medskip

\noindent{\bf Remarks}
\begin{itemize}
\item The case $f$ diagonal but  singular can be treated in the same way. In fact $Ker(f)$ is an abelian ideal of $\g$.
\item If $f$ is diagonalizable, up to an automorphism of $\h_{2p+1}$ we can return to the previous case.
\end{itemize}

\section{Classification of contact Lie algebras in small dimensions}

\subsection{Classification of real contact $3$-dimensional Lie algebra}

Assume that $\g$ is a linear deformation of $\h_3$. Then it is a central extension of a symplectic $2$-dimensional Lie algebra.   We deduce that $\varphi_1=0$ or $\varphi_1(X_1,X_2)=X_1.$ In the first case $\g$ is isomorphic to $\h_3$. In the second case the Lie bracket of $\g$ is isomorphic to
$$
\left\{
\begin{array}{l}
\medskip
[X_{1},X_{2}]=X_{3}+X_1,\\
\medskip
\lbrack X_2,X_3 \rbrack=0.
\end{array}
\right.
$$
\medskip

Assume that $\g$ is a quadratic deformation of $\h_3$. Its corresponds an endomorphism $f$ of $\frak{k}$ generated by $\{X_1,X_2\}$ whose matrix is
$$
\left(
       \begin{array}{cc}
         a & b \\
         c & -a \\
       \end{array}
     \right)
$$
and $\varphi_2$ is given by
$$
\left\{
\begin{array}{l}
\varphi_2(X_1,X_3)=aX_1+cX_2,\\
\varphi_2(X_2,X_3)=bX_1-aX_2,\\
\end{array}
\right.
$$
and $\delta_{\mu_0}\varphi_2=\varphi_2 \circ \varphi_2=0.$ If $\varphi_1(X_1,X_2)=\alpha X_1+\beta X_2$, then
the identity $\varphi_1\circ\varphi_2 +\varphi_2\circ\varphi_1=0$ is equivalent to
$$\alpha f(X_1)+\beta f(X_2)=0.$$

If $\varphi_1=0$, we have no conditions on $f$. If $rank(f)=2$, thus $a^2+bc \neq 0$, the eigenvalues of $f$ are $\lambda_1$ and $-\lambda_1$ and $f$ is diagonalizable. Let $\sigma$ be an automorphism of $\h_3$. Its matrix in the basis $\{X_1,X_2,X_3\}$ is
$$\left(
    \begin{array}{ccc}
      a_{11} & a_{12} & 0 \\
      a_{21} & a_{22} & 0 \\
      a_{31} & a_{32} & a_{11}a_{22} -a_{21}a_{12}\\
    \end{array}
  \right)
$$
Its restriction to the space $\frak{k}$ generated by $X_1$ and $X_2$ is a non singular matrix and conversely any non singular matrix of order $2$ defines an automorphism of $\h_3$.  Then we can reduce the matrix of $f$ in the diagonal form modulo an automorphism of $\h_3$ and assume that the matrix of $f$ is in the diagonal form
$$
f=\left(
       \begin{array}{cc}
         \lambda_1 & 0 \\
         0 & -\lambda_1\\
       \end{array}
     \right)
$$
with $\lambda_1\neq 0.$ If $\lambda_1 \in \R$, the Lie bracket of $\g$ is
$$
[X_1,X_2]=X_3, \ [X_1,X_3]=\lambda_1X_1, \ [X_2,X_3]=-\lambda_1X_2$$
and $\g$ is isomorphic to $sl(2,\R)$. If $\lambda_1 $ is a complex number, then  $\g$ is
 isomorphic to $so(3).$ Then any simple real Lie algebra of rank $1$ admits a contact form.

If $\varphi_1 \neq 0$, we can take, up to an automorphism of $\h_3$, $\varphi_1(X_1,X_2)=X_1$. This implies that 
$$(\varphi_1 \circ \varphi_2+\varphi_2 \circ \varphi_1)(X_1,X_2,X_3)=f(X_1)=0$$
 and the matrix of $f$ is
$$
\left(
       \begin{array}{cc}
         0& b \\
         0 & 0 \\
       \end{array}
     \right)
$$
with $b \neq 0$. We deduce that  $\varphi_2(X_1,X_3)=0, \ \varphi_2(X_2,X_3)=bX_1$ and the Lie bracket of $\g$ is
$$
\left\{
\begin{array}{l}
\medskip
[X_{1},X_{2}]=X_{3}+X_1,\\
\medskip
\lbrack X_2,X_3 \rbrack=bX_1.
\end{array}
\right.
$$

\subsection{Classification in dimension $5$.}
Assume that $\g$ is a quadratic deformation of $\h_5$. 
The matrix of $f$ is
$$\left(
    \begin{array}{cccc}
      a_1 & b_1 & -b_4 & b_3 \\
      a_2 & -a_1 & a_4 & -a_3 \\
      a_3 & b_3 & c_3 & d_3 \\
      a_4 & b_4 & c_4 & -c_3 \\
    \end{array}
  \right)
  $$
  \begin{description}
    \item[1. $f$ is diagonalizable]
  \end{description}
  \noindent Its eigenvalues are $\lambda_1,-\lambda_1,\lambda_2,-\lambda_2.$ There exists a change of basis given by an automorphism of $\h_5$ which reduces the matrix of $f$ in a diagonal form. In fact, for any square matrix, we denote by $\Delta_{ij}^{kl}$ the determinant of the submatrix constituted of elements belonging to the lines $i,j$ and columns $k,l$. Any automorphism of $\h_5$ is given by a matrix, related with the classical basis, such that
  $$
  \begin{array}{l}
    \Delta_{12}^{12} +\Delta_{34}^{12}= \Delta_{12}^{34} +\Delta_{34}^{34} \neq 0, \\
    \Delta_{12}^{13} +\Delta_{34}^{13} = \Delta_{12}^{14} +\Delta_{34}^{14}=\Delta_{12}^{23} +\Delta_{34}^{23}=\Delta_{12}^{24} +\Delta_{34}^{24}=0.
  \end{array}
  $$
  For any automorphism $\sigma$, let $\widetilde{\sigma}$ be the isomorphism of $\R^4$ whose matrix is the matrix of $\sigma$ where we have removed the line $5$ and the column $5$. The previous conditions are still valide for  $\widetilde{\sigma}$. Thus, we can find a new classical basis of $\h_5$ such that $f(X_1)=\lambda_1X_1,f(X_2)=-\lambda_1X_2.$ The matrix of $f$ is now of the form
  $$\left(
    \begin{array}{cccc}
      \lambda_1 & 0 & 0 & 0 \\
      0 & -\lambda_1 & 0 & 0 \\
      0 & 0 & c_3 & d_3 \\
      0 & 0 & c_4 & -c_3 \\
    \end{array}
  \right)
  $$
  and using the same arguments, we can reduce this matrix in the diagonal form in a classical basis of $\h_5$, always denoted $\{X_1,X_2,X_3,X_4,X_5\}$, that is,
  $$f(X_1)=\lambda_1X_1,f(X_2)=-\lambda_1X_2,f(X_3)=\lambda_2X_3,f(X_4)=-\lambda_2X_4.$$
  If we write $\varphi_1(X_i,X_j)=\sum_{k=1}^4 A_{ij}^kX_k$, thus the equations $\varphi_1\circ\varphi_2+\varphi_2\circ\varphi_1(X_i,X_j,X_5)=0$ for $1 \leq i < j \leq 4$ give
  $$
\left\{
       \begin{array}{l}
       \lambda_1A_{12}^1=\lambda_1A_{12}^2=\lambda_2A_{12}^3=\lambda_2A_{12}^4=0,\\
        -\lambda_2A_{13}^1=-(2\lambda_1+\lambda_2)A_{13}^2=-\lambda_1A_{13}^3=-(\lambda_1+2\lambda_2)A_{13}^4=0,\\
         \lambda_2A_{14}^1=(-2\lambda_1+\lambda_2)A_{14}^2=(-\lambda_1+2\lambda_2)A_{14}^3=-\lambda_1A_{14}^4=0,\\
         (2\lambda_1-\lambda_2)A_{23}^1=-\lambda_2A_{23}^2=\lambda_1A_{23}^3=(\lambda_1-2\lambda_2)A_{14}^4=0,\\
         (2\lambda_1+\lambda_2)A_{24}^1=\lambda_2A_{24}^2=(\lambda_1+2\lambda_2)A_{24}^3=\lambda_1A_{24}^4=0,\\
         \lambda_1A_{34}^1=-\lambda_1A_{34}^2=\lambda_2A_{34}^3=-\lambda_2A_{34}^4=0.
           \end{array}
  \right.
  $$
  i) If $\lambda_1\neq 0, \ \lambda_2\neq 0, \ \lambda_2\neq \pm 2\lambda_1, \ \lambda_1\neq \pm 2\lambda_2$, then $\varphi_1=0$.

  \noindent ii) If $\lambda_1= 0, \ \lambda_2\neq 0$, then the equations $\varphi_1\circ\varphi_1 +\delta\varphi_2=0$ and $\delta\varphi_1=0$ give the following Lie algebras
  $$
\left\{
\begin{array}{l}
\medskip
[X_{1},X_{2}]=X_{5}+aX_1+bX_2,\
\lbrack X_1,X_3 \rbrack=cX_3, \ \lbrack X_1,X_4 \rbrack=-cX_4,\\
\medskip
\lbrack X_2,X_3 \rbrack=dX_3, \ \lbrack X_2,X_4 \rbrack=-dX_4, \ \lbrack X_3,X_4 \rbrack=X_5,\\
\medskip
\lbrack X_3,X_5 \rbrack=(ac+bd)X_3, \ \lbrack X_4,X_5 \rbrack=-(ac+bd)X_4,\\
\end{array}
\right.
$$
with $a,b,c, d \in \R$;
 $$
\left\{
\begin{array}{l}
\medskip
[X_{1},X_{2}]=X_{5}+bX_2, \
\lbrack X_1,X_3 \rbrack=cX_3, \ \lbrack X_1,X_4 \rbrack=(b-c)X_4,\\
\medskip
\lbrack X_2,X_3 \rbrack=dX_3, \ \lbrack X_2,X_4 \rbrack=-dX_4, \ \lbrack X_3,X_4 \rbrack=bX_2+X_5,\\
\medskip
\lbrack X_3,X_5 \rbrack=bdX_3, \ \lbrack X_4,X_5 \rbrack=-bdX_4,\\
\end{array}
\right.
$$
with $b,c, d \in \R;$
$$
\left\{
\begin{array}{l}
\medskip
[X_{1},X_{2}]=X_{5}+aX_1+bX_2, \
\lbrack X_1,X_3 \rbrack=cX_3, \ \lbrack X_1,X_4 \rbrack=(b-c)X_4,\\
\medskip
\lbrack X_2,X_3 \rbrack=dX_3, \ \lbrack X_2,X_4 \rbrack=(-a-d)X_4, \ \lbrack X_3,X_4 \rbrack=aX_1+bX_2+X_5,\\
\medskip
\lbrack X_3,X_5 \rbrack=(ac+bd)X_3, \ \lbrack X_4,X_5 \rbrack=-(ac+bd)X_4,\\
\end{array}
\right.
$$
with $b,c, d \in \R.$

\noindent iii) If $\lambda_2= 0, \ \lambda_1\neq 0$, we find a case is isomorphic to the case ii).

\noindent iv) If $\lambda_2=2\lambda_1$, we have $\varphi_1(X_1,X_4)=A_{14}^2X_2, \varphi_1(X_2,X_3)=A_{23}^1X_1$ and $\varphi_1(X_i,X_j)=0$ in all the other cases. We deduce
$$\varphi_1\circ\varphi_1(X_1,X_2,X_3)=0=-\delta\varphi_2(X_1,X_2,X_3)=-\lambda_2X_3.$$
Thus $\lambda_2=0$, which is impossible. We have similar results for $\lambda_2=-2\lambda_1$ and $\lambda_1=\pm 2\lambda_2$.

 \begin{description}
    \item[2. $f$ is not diagonalizable]
  \end{description}
Since the eigenvalues of $f$ are of the form $\lambda_1,-\lambda_1,\lambda_2,-\lambda_2$, if one of these eigenvalues is complex, for example $\lambda_1$, thus $\overline{\lambda_1}=-\lambda_1$ and $\lambda_1=ia$ with $a \in \R$. In a first step assume that $\lambda_1 = ia$ and $\lambda_2 \in \R$. Up to the restriction to $\R^4$ of an automorphism of $\h_5$, we can assume that the matrix of $f$ is
 $$\left(
    \begin{array}{cccc}
     0 & a & 0 & 0 \\
      -a& 0 & 0 & 0 \\
      0 & 0 & \lambda_2 &0\\
      0 & 0 & 0 & -\lambda_2 \\
    \end{array}
  \right)
  $$
 with $a, \lambda_2 \in \R$.

 \noindent i) If $a \neq 0$ and $\lambda_2 \neq 0$, then $\varphi_1=0.$

  \noindent ii) If $a \neq 0$ and $\lambda_2 = 0$,
 the relation $\varphi_1\circ\varphi_1 +\delta\varphi_2=0$ gives
  $$
\left\{
\begin{array}{l}
\varphi_1(X_1,X_2)=A_{12}^3X_3+A_{12}^4X_4,\\
\varphi_1(X_1,X_3)=A_{13}^1X_1+A_{13}^2X_2,\\
\varphi_1(X_2,X_3)=-A_{13}^2X_1+A_{13}^1X_2,\\
\varphi_1(X_1,X_4)=A_{14}^1X_1+A_{14}^2X_2,\\
\varphi_1(X_2,X_4)=-A_{14}^2X_1+A_{14}^1X_2,\\
\varphi_1(X_3,X_4)=A_{34}^1X_1+A_{34}^2X_2.\\
\end{array}
  \right.
  $$
  The relation $\delta\varphi_1=0$ implies $A_{12}^3=-2A_{14}^1,A_{12}^4=-2A_{13}^1$ and $(\varphi_1\circ\varphi_1 +\delta\varphi_2)(X_1,X_2,X_3)=0$ is equivalent to
   $$
\left\{
\begin{array}{l}
2 A_{13}^1A_{34}^3-2 A_{13}^1A_{14}^1- 2A_{13}^1A_{14}^1+\lambda_2=0,\\
2 A_{13}^1A_{34}^4-4 A_{13}^1A_{13}^1=0,
\end{array}
  \right.
  $$
  and $(\varphi_1\circ\varphi_1 +\delta\varphi_2)(X_1,X_2,X_4)=0$ induces that
    $$
\left\{
\begin{array}{l}
2 A_{14}^1A_{34}^3-4 A_{14}^1A_{14}^1=0,\\
2 A_{14}^1A_{34}^4+2 A_{14}^1A_{13}^1+2 A_{14}^1A_{13}^1-\lambda_2=0.
\end{array}
  \right.
  $$
  \begin{itemize}
    \item $A_{14}^1=A_{13}^1=0$. We obtain the following Lie algebras
    $$
\left\{
\begin{array}{l}
\medskip
[X_{1},X_{2}]=X_{5}, \lbrack X_3,X_4 \rbrack=eX_3+fX_4+X_5,\\
\medskip
\lbrack X_1,X_3 \rbrack=cX_2, \ \lbrack X_1,X_4 \rbrack=dX_2,\\
\medskip
\lbrack X_2,X_3 \rbrack=-cX_1, \ \lbrack X_2,X_4 \rbrack=-dX_1,\\
\medskip
\lbrack X_1,X_5 \rbrack=-(ce+df)X_2, \ \lbrack X_2,X_5 \rbrack=(ce+df)X_1.\\
\end{array}
\right.
$$
    \item $A_{13}^1=A_{34}^4=0, \ A_{34}^3=2A_{14}^1$. We obtain the following Lie algebras
    $$
\left\{
\begin{array}{l}
\medskip
[X_{1},X_{2}]=2aX_3+X_{5}, \lbrack X_3,X_4 \rbrack=2aX_3+X_5,\\
\medskip
\lbrack X_1,X_3 \rbrack=cX_2, \ \lbrack X_1,X_4 \rbrack=aX_1+dX_2,\\
\medskip
\lbrack X_2,X_3 \rbrack=-cX_1, \ \lbrack X_2,X_4 \rbrack=-dX_1+aX_2,\\
\medskip
\lbrack X_1,X_5 \rbrack=-2acX_2, \ \lbrack X_2,X_5 \rbrack=2acX_1.\\
\end{array}
\right.
$$
    \item $A_{14}^1=A_{34}^3=0, \ A_{34}^4=-2A_{13}^1$. This case is isomorphic to the previous case.
    \item $ A_{34}^3=2A_{14}^1, \ A_{34}^4=-2A_{13}^1$. We obtain the following Lie algebras
    $$
\left\{
\begin{array}{l}
\medskip
[X_{1},X_{2}]=2aX_3-2bX_4+X_{5}, \lbrack X_3,X_4 \rbrack=2aX_3-2bX_4+X_5,\\
\medskip
\lbrack X_1,X_3 \rbrack=bX_1+cX_2, \ \lbrack X_1,X_4 \rbrack=aX_1+dX_2,\\
\medskip
\lbrack X_2,X_3 \rbrack=-cX_1+bX_2, \ \lbrack X_2,X_4 \rbrack=-dX_1+aX_2,\\
\medskip
\lbrack X_1,X_5 \rbrack=-2(ac-bd)X_2, \ \lbrack X_2,X_5 \rbrack=2(ac-bd)X_1.\\
\end{array}
\right.
$$
  \end{itemize}

\noindent iii) If $a = 0$ and $\lambda_2 \neq 0$, $f$ is diagonal.

\noindent iv) Assume now that $\lambda_1$ is a real eigenvalue of $f$ of multiplicity $2$.  If $\{X_1,X_2,X_3,X_4\}$ is a Jordan basis, that is, 
$$f(X_1)=\lambda_1X_1,\ f(X_2)=-\lambda_1X_2, \   f(X_3)=\lambda_1X_3+X_1, \  f(X_4)=-\lambda_1X_4+X_2, $$
then the equation $(\varphi_1 \circ \varphi_2 +\varphi_2\circ  \varphi_1)(X_i,X_j,X_5)=0$ implies $\varphi_1=0$  if $\lambda_1 \neq 0$. We have already studied this case. If $\lambda_1=0$, then 
$$\varphi_2(X_3,X_5)=X_1, \ \varphi_2(X_4,X_5)=X_2,  \varphi_2(X_1,X_5)=\varphi_2(X_2,X_5)=0.$$
We deduce $\delta_{\mu_0}\varphi_2(X_1,X_4,X_5)=X_5$ and we cannot have $\varphi_1 \circ \varphi_1+\delta_{\mu_0}\varphi_2=0.$ Then this case is impossible.

\section{Contact nilpotent Lie algebras}
Recall that we call contact Lie algebra, an odd dimensional Lie algebra $\g$ provided with a linear form whose Cartan class is equal to the dimension of this Lie algebra. This form is also called a contact form on $\g$.

\subsection{Contact $2$-step nilpotent  Lie algebras}
A nilpotent Lie algebra is $2$-step nilpotent if its nilindex is equal to $2$ and the bracket $\mu$ satisfies $\mu(\mu(X,Y),Z)=0$. In particular the Heisenberg algebras are $2$-step nilpotent. Let $\g$ be a $(2p+1)$-dimensional Lie algebra with a contact form. Its bracket $\mu$ is a deformation of the bracket $\mu_0$ of the Heisenberg algebra $\h_{2p+1}$. But, from \cite{GR-Lieass}, $\mu_0$ is rigid in the variety $2\mathcal{N}ilp(2p+1)$ of $(2p+1)$-dimensional $2$-step nilpotent Lie algebras. This means that any deformation of $\mu_0$ in a  $2$-step nilpotent Lie algebra is isomorphic to $\mu_0$. Then $\g$ is isomorphic to $\h_{2p+1}$.
\begin{proposition}
Any $(2p+1)$-dimensional contact $2$-step nilpotent Lie algebra  is isomorphic to $\h_{2p+1}$.
\end{proposition}

\subsection{Contact $p$-step nilpotent  Lie algebras}

\begin{lemma}\label{center}\cite{GozeCras1}
The center $Z(\g)$ of a contact nilpotent Lie algebra $\g$ is of dimension $1$.
\end{lemma}
\pf Since $\g$ is nilpotent, its center $Z(\g)$ is not trivial. Let $\omega \in \g^*$ be a contact form on $\g$. Any vector $X \in Z(\g)$,  $X \neq 0$ satisfies $i(X)(d\omega)=0$. In fact 
$$i(X)d\omega(Y)=-\omega [X,Y]=0$$
for any $Y \in \g$. Then the vector space
$$A(d\omega)=\{X \in \g, \ i(X)d\omega=0\}$$
contains $Z(\g)$. Since $\omega$ is a contact form, the characteristic space $\mathcal{C}(\omega)$ which is equal by definition to $A(\omega)\cap A(d\omega)$ is equal to ${0}$. This implies $\dim A(d\omega) \leq 1$. Then $A(d\omega)=Z(\g)$ and is one-dimensional.

\begin{theorem}\label{contactnilpotent}
Let $\g$ be a $(2p+1)$-dimensional $k$-step nilpotent Lie algebra. Then there exists on $\g$ a contact form if and only if $\g$ is a central extension of a $(2p)$-dimensional $(k-1)$-step  nilpotent symplectic Lie algebra $\h$, the extension being defined by the $2$-cocycle  given by the symplectic form.
\end{theorem}
\pf Let $\g$ be a $(2p+1)$-dimensional nilpotent Lie algebra and $\omega$ a contact form on $\g$. From the Lemma \ref{center}, $A(d\omega)=Z(\g)$ is a one-dimensional space so $\g= Ker(\omega) \oplus Z(\g)$. Thus we can find a basis $\{X_1,\cdots,X_{2p+1}\}$ associated with this decomposition and corresponding to a basis satisfying (\ref{basis}).  The Lie algebra $\g$ is isomorphic to a quadratic deformation of $\h_{2p+1}$
$$\mu=\mu_0+t\varphi_1+t^2\varphi_2$$
the bilinear forms $\varphi_1$ and $\varphi_2$ satisfying (\ref{phi1}) and (\ref{phi2}). In particular 
$$\varphi_2(X_i,X_{2p+1})=\mu(X_i,X_{2p+1})=0.$$This implies $\varphi_2=0$ and any contact nilpotent Lie algebra is isomorphic to a linear deformation of the Heisenberg algebra. Now we assume that $\mu=\mu_0+t\varphi_1$ is $k$-step nilpotent. If $\psi_1$ and $\psi_2$ are bilinear map on $\g$, we denote by $\psi_1 \bullet \psi_2$ the trilinear map
$$\psi_1 \bullet \psi_2(X,Y,Z)=\psi_1(\psi_2(X,Y),Z)$$
for any $X,Y,Z \in \g$. This product is sometimes denoted by $\circ_1$ in the Gerstenhaber terminology. Since $\mu$ is $k$-step nilpotent, we have $0=\mu^{\bullet^k}=\mu \bullet \cdots \bullet \mu$ ($k$-times). We develop the identity $(\mu_0+t\varphi_1)^{\bullet^k}=0$ and we reduce the result by the relations $\mu_0 \bullet \mu_0=0$ because $\mu_0$ is $2$-step nilpotent and $\varphi_1 \bullet \mu_0=0$ from (\ref{phi1}). We obtain $\mu_0 \bullet \varphi_1^{\bullet^{k-1}}=0$ and $ \varphi_1^{\bullet^{k}}=0.$ But $\mu_0 \bullet \varphi_1^{\bullet^{k-1}}=0$ is equivalent to
$$\mu_0 ( \varphi_1^{\bullet^{k-1}}(Y_1,\cdots,Y_{k-1}),Y)=0$$
for all $Y$. This implies $\varphi_1^{\bullet^{k-1}}(Y_1,\cdots,Y_{k-1}) \in Z(\h_{2p+1})$. From (\ref{phi1}), 
$$\varphi_1^{\bullet^{k-1}}(Y_1,\cdots,Y_{k-1})=0.$$
But $\varphi_1$ is also a Lie bracket. Then it defines a $(2p)$-dimensional $(k-1)$-step nilpotent Lie algebra. The result follows from  Theorem \ref{symplec}.

\subsection{Contact filiform Lie algebras.}


Let $\g$ be a $(2p+1)$-dimensional filiform Lie algebra.  It is equivalent to say that $\g$ is $(2p)$-step nilptent. From Theorem \ref{contactnilpotent}, $\g$ is a central extension of a $(2p)$-dimensional $(2p-1)$-step  nilpotent symplectic Lie algebra $\h$, that is, $\g$ is a central extension of a symplectic filiform Lie algebra of dimension $2p$.

In \cite{GozeKhakimMedina}, we have done a direct construction of contact filiform Lie algebras. We recall quickly this construction. Let $\g$ be a $(2p+1)$-dimensional filiform Lie algebra.  Let $\{e_0, \cdots, e_{2p} \}$ an adapted basis of $\g$ (see \cite{GozeKhakimbook}). This basis satisfies in particular $[e_0, e_i]=e_{i+1}, i=1,\cdots,2p-1$ and $e_{2p}$ is a basis of $Z(\g)$. If $\omega$ is a contact form, then from Lemma \ref{center}, $\omega(e_{2p})\neq 0.$  If $\{ \omega_0, \cdots ,\omega_{2p}\}$ is the dual basis of $\{ e_0, \cdots , e_{2p} \}$ and if $\g$ is a contact Lie algebra, then $\omega_{2p}$ is also a contact form (see \cite{GozeKhakimMedina}).
 Let us consider the following filiform Lie algebra whose bracket $\mu$ writes
 $$ \mu=\mu_{L_{2p+1}}+a_{1,4}\psi_{1,4} + \cdots +a_{i,2i+2}\psi_{i,2i+2} +\cdots +a_{p-1,2p}\psi_{p-1,2p}
$$
with $a_{i,2i+2}\neq 0$ for $i=1, \cdots ,p$
and where $\mu_{L_{2p+1}}$ is the Lie bracket of $L_{2p+1}$, that is,
$$ \mu_{L_{2p+1}}(e_0,e_i)=e_{i+1} \ {\rm for \ } i=1, \cdots, 2p-1$$
other brackets areequal to $0$. The bilinear forms $\psi_{k,s}$ are given by
$$
\left\{
\begin{array}{l}
\psi_{k,s}(e_i,e_j)=(-1)^k C^{j-k-1}_{k-i} (ad \ e_0)^{i+j-2k-1}(e_s), \ \ 1 \leq i < j-1 \leq n-1,\\ 
\psi_{k,s}(e_k,e_j)=e_{s+j-k-1}
\end{array}
\right.
$$
 and other non given values of $\psi_{k,s}$ are zero (or defined by skewsymmetry).
Thus $\mu$ is the Lie bracket of a filiform contact Lie algebra if and only if $A_1 \cdots A_{p-1}\neq 0$ with
$$A_i=\sum_{k=0}^{i-1} (-1)^k a_{p-1+k, 2p-2i+2k+2} C^{k}_{2i-k-2}.$$
For example, if $2p+1=9$, we have 
$$\mu_{c,9}=\mu_{L_9}+a_{1,4}\psi_{1,4} +a_{2,6}\psi_{2,6} +a_{3,8}\psi_{3,8}$$
with 
$$\left\{ 
\begin{array}{l}
a_{3,8} \neq 0, \\
a_{2,6} - a_{3,8} \neq 0, \\
a_{1,4} - 3 a_{2,6} + a_{3,8} \neq 0.  
\end{array}
\right.
$$
More precisely 
$$\left\{
\begin{array}{l}
\mu_{c,9}(e_0,e_i)=e_{i+1}, \ i=1, \cdots ,7, \\
\mu_{c,9}(e_1,e_2)= a_{1,4}e_4,\\
\mu_{c,9}(e_1,e_3)= a_{1,4}e_5, \\
\mu_{c,9}(e_1,e_4)= (a_{1,4}-a_{2,6})e_6, \\
\mu_{c,9}(e_1,e_5)= (a_{1,4}-2a_{2,6})e_7, \\
\mu_{c,9}(e_1,e_6)= (a_{1,4} - 3 a_{2,6} + a_{3,8})e_8, \\
\mu_{c,9}(e_2,e_3)= a_{2,6}e_6, \\
\mu_{c,9}(e_2,e_4)= a_{2,6}e_7,\\
\mu_{c,9}(e_2,e_5)= (a_{2,6}-a_{3,8})e_8, \\
\mu_{c,9}(e_3,e_4)= a_{3,8}e_8.
\end{array}
\right.
$$

\begin{proposition}
Any $(2p+1)$-dimensional filiform contact Lie algebra is a linear deformation of the filiform Lie algebra $\mathfrak{k}$ whose bracket is isomorphic to $\mu_{c,2p+1}.$ 
\end{proposition}

Let us remark that the Lie algebra $\mathfrak{k}$ is characteristically nilpotent but admits a $\mathbb{Z}_2$-grading \cite{BGR}.

\subsection{Symplectic Lie algebras}

Theorem \ref{contactnilpotent} shows that the construction of contact nilpotent Lie algebras is equivalent to the construction of symplectic nilpotent Lie algebras. In \cite{Medina} a process called double extension is described to construct all the $(2p)$-dimensional symplectic Lie algebra from the class of $(2p-2)$-dimensional symplectic  Lie algebras.
Theorem  \ref{contactnilpotent} gives an other characterization of symplectic Lie algebras.
 Let $\varphi$ be a $2$-cocycle of $\h_{2p+1}$ satisfying $\varphi(X,X_{2p+1}) =0$ for any $X$ and $\varphi(X,Y)$ is in the 
 subspace $E$ generated by $\{X_1, \cdots, X_{2p}\}.$ This cocycle defines a bilinear map $\psi$ on $E$ and $\psi$ satisfies the Jacobi conditions, that is, $\psi \circ \psi=0$ then the Lie algebra $(E, \psi)$ is a $(2p)$-dimensional symplectic Lie algebra. Conversely if $\frak{m}$ is a $(2p)$-dimensional symplectic Lie algebra such that the symplectic form $\theta$ is not exact, then $\theta$ defines a central extension of $\frak{m}$ and the obtained Lie algebra $\g$ admits a contact form.
 Moreover $\g$ is isomorphic to a linear deformation of $\h_{2p+1}$ and it is given by a $2$-cocycle of $\h_{2p+1}$.
 
 \begin{theorem}
 The Lie bracket of a $(2p)$-dimensional symplectic Lie algebra is given by the restriction of a $2$-cocycle of $\h_{2p+1}$ to the $(2p)$-dimensional subspace $E$ of $\h_{2p+1}$ whose kernel contains the center of $\h_{2p+1}$.
 
 Conversely  if such restriction satisfies Jacobi identity then it is the Lie bracket of a symplectic Lie algebra.
 \end{theorem}

In the following section we will describe symplectic Lie algebras with an exact symplectic form.

\section{Frobeniusian Lie algebras}
A $2p$-dimensional Lie algebra is called frobeniusian if there exists $\omega \neq 0, \omega \in \g^*$ such that $(d\omega)^p \neq 0.$  We have seen that nilpotent and semisimple Lie algebras are never frobeniusian. In this section, we will describe $2p$-dimensional frobeniusian Lie algebras which give, by deformation, any fobeniusian Lie algebras.

\subsection{Contractions of Lie algebras}
Let $\g_0$ be a $n$-dimensional Lie algebra  whose Lie bracket is denoted by $\mu_0$. We consider $\left\{  f_{t}\right\} _{t \in ]0,1]} $ a sequence of isomorphisms in
$\mathbb{K}^{n}$ with $\mathbb{K}=\R$ or $\C$. Any Lie bracket
\[
\mu_{t}=f_t^{-1}\circ\mu_{0}\left(  f_{t}\times f_{t}\right)
\]
corresponds to a Lie algebra $\g_t$ which is isomorphic to $\g_0$. If the limit
$\lim_{t\rightarrow 0}\mu_{t}$ exists (this limit is computed in the finite dimensional vector space of bilinear maps in $\K^n$), it defines a Lie bracket $\mu$ of a $n$-dimensional Lie algebra $\g$
called a contraction of $\g_{0}$.

\noindent{\bf Remark}. Let $\frak{L}^{n}$ be the variety of Lie algebra laws over $\mathbb{C}^{n}$ provided with its Zariski topology. The algebraic structure of this variety is defined by the Jacobi polynomial equations on the structure constants.
The linear group $GL\left(  n,\mathbb{C}\right)  $ acts on
$\mathbb{C}^{n}$ by changes of basis. A Lie algebra $\frak{g}$ 
is contracted to the $\g_0$ if the Lie bracket of $\g$ 
 is in the closure of the orbit of the Lie bracket of $\g_ 0$ by the group
action (for more details see \cite{GozeEllipse}).

\begin{definition}\cite{GozeCras2}
The algebra $\frak{g}_{0}$ is called a semimodel relative to property $\left(
P\right)  $ if any Lie algebra law $\mu$ satisfying $\left(  P\right)  $
contracts on $\mu_{0}$.\newline The semimodel is called a model if, in
addition, any deformation of $\mu_{0}$ satisfies $\left(  P\right)  $. We use the same vocabulary even if instead of an algebra $\g_0$ we have a family of algebras irreducible by contractions.
\end{definition}
This means that the neighbourhoods in $\frak{L}^{n}$ of a model related to a given property $(P)$  are entirely characterized by property
$\left(  P\right)  $. For example, if $(P)$ is the property "{\it there exists a contact form}", then, from the previous section, we obtain for any odd dimension a unique model, the Heisenberg algebra.

\subsection{Classification of frobeniusian Lie algebras up to a contraction}
Let $\g$ be a $2p$-dimensional frobeniusian Lie algebra. There exists a basis $\{X_1, \cdots, X_{2p}\}$ of $\g$ such that the dual basis $\{\omega_1,\cdots,\omega_{2p}\}$ is adapted to the frobeniusian structure, that is, $$d\omega_1=\omega_1\wedge \omega_2+\omega_3\wedge \omega_4+\cdots+\omega_{2p-1}\wedge \omega_{2p}.$$ In the following, we define Lie algebras, not with its brackets, but with its Maurer-Cartan equations. We assume here that $\K=\C$.\begin{theorem}\label{frob}
Let $\g_{a_1,\cdots,a_{p-1}}$, \ $ a_i \in \C$ be    the Lie algebras defined by
\[
\left\{
\begin{array}
[c]{l}%
d\omega_{1}=\omega_{1}\wedge\omega_{2}+\sum_{k=1}^{p-1}\omega_{2k+1}%
\wedge\omega_{2k+2},\\
d\omega_{2}=0,\\
d\omega_{2k+1}=a_{k}\omega_{2}\wedge\omega_{2k+1},\;1\leq k\leq p-1,\\
d\omega_{2k+2}=-\left(  1+a_{k}\right)  \omega_{2}\wedge\omega
_{2k+2},\;1\leq k\leq p-1.
\end{array}
\right.
\]
Then any $2p$-dimensional frobeniusian Lie algebra is contracted in an element of the family $\mathcal{F}=\{\g_{a_1,\cdots,a_{p-1}}\}_{a_i \in \C}$. Moreover, any element of $\mathcal{F}$ cannot be contracted in other different element of this family.
\end{theorem}
This means that the family $\mathcal{F}$ is a model for the property to be frobeniusian.
\pf It can be easily seen that the Jacobi conditions are satisfied for any $\{a_1,\cdots,a_k\}$. Let $\g$ be a frobeniusian Lie algebra. There exists a basis of $\g^*$ such that
$$d\omega_{1}=\omega_{1}\wedge\omega_{2}+\sum_{k=1}^{p-1}\omega_{2k+1}\wedge\omega_{2k+2},$$
and the frobeniusian structure is given by $\omega=\omega_1.$ Let $\{X_1,\cdots,X_{2p}\}$ the dual basis. We consider the linear isomorphisms $f_\varepsilon$ with $f_\varepsilon(X_1)=\varepsilon^2 X_1,$ $f_\varepsilon(X_i)=\varepsilon X_i$ for $i \neq 1,2$ and $f_\varepsilon(X_2)=X_2$. In the Lie algebra $\g_\varepsilon$ we have
$$d\omega_{1}=\omega_{1}\wedge\omega_{2}+\sum_{k=1}^{p-1}\omega_{2k+1}\wedge\omega_{2k+2}$$
and if $C_{ij}^k$ are the structure constants of $\g$:
\[
\left\{
\begin{array}{l}
d\omega_{2}=\varepsilon^2 C_{12}^2\omega_{1}\wedge\omega_{2}+\varepsilon\sum C_{2,i}^2\omega_{2}\wedge\omega_{i}+\varepsilon^2\sum C_{ij}^2\omega_{i}\wedge\omega_{j},\\
d\omega_{2k+1}=\sum_{i \neq 1} C_{2,i}^{2k+1}\omega_{2}\wedge\omega_{i}
+\varepsilon(C_{12}^{2k+1}\omega_{1}\wedge\omega_{2}+\sum_{2 <i<j<2p} C_{ij}^{2k+1}\omega_{i}\wedge\omega_{j},\\
d\omega_{2k+2}=\sum_{i \neq 1}C_{2,i}^{2k+2}\omega_{2}\wedge\omega_{i}
+\varepsilon(C_{12}^{2k+2}\omega_{1}\wedge\omega_{2}+\sum_{2 <i<j<2p} C_{ij}^{2k+2}\omega_{i}\wedge\omega_{j}).
\end{array}
\right.
\]
When $\varepsilon \rightarrow 0$, we obtain the Lie algebras whose  Maurer Cartan equations are
\begin{equation}\label{7}
\left\{
\begin{array}{l}
d\omega_{1}=\omega_{1}\wedge\omega_{2}+\sum_{k=1}^{p-1}\omega_{2k+1}\wedge\omega_{2k+2},\\
d\omega_{2}=0,\\
d\omega_{2k+1}=\sum_{i \neq 1} C_{2,i}^{2k+1}\omega_{2}\wedge\omega_{i}
,\\
d\omega_{2k+2}=\sum_{i \neq 1}C_{2,i}^{2k+2}\omega_{2}\wedge\omega_{i}.
\end{array}
\right.
\end{equation}
The Jacobi conditions writes
$$
\left\{
\begin{array}{l}
C_{2,2s}^{2s}=-1-C_{2,2s-1}^{2s-1}\\
C^{2i+1}_{2,2k+1}=-C^{2k+2}_{2,2i+2}, \  C^{2i+1}_{2,2k+2}=C^{2k+1}_{2,2i+2},\\
C^{2i+2}_{2,2k+1}=C^{2k+2}_{2,2i+1}, \  C^{2i+1}_{2,2k+2}=-C^{2k+1}_{2,2i+1}.
\end{array}
\right.
$$
It remains to show that Lie algebras (\ref{7}) can be contracted on the models described in the Theorem.   To simplify notations we put $a_i^j=C^j_{2,i}$ and $M=(a_i^j)$, $3 \leq i,j \leq 2p$. Let $V$ be the vector space generated by $\{X_3,\cdots,X_{2p}\} $,  and $f$ the morphism of $V$ whose matrix in the given basis is $M$. If we denote by $[X,Y]$ the bracket of any Lie algebras (\ref{7}), then $[V,V] \subset \C{X_1}.$ We consider the following steps:
\begin{itemize}
\item Step 1: Complex classification of the Lie algebras (\ref{7}). These Lie algebras are in correspondence with the matrix $M$. Then we have to reduce $M$ modulo the subgroup of $GL(2p,\C)$ which is invariant by the frobeniusian form $d\omega_{1}=\omega_{1}\wedge\omega_{2}+\sum_{k=1}^{p-1}\omega_{2k+1}\wedge\omega_{2k+2}.$ 
\begin{itemize}
\item Let $\lambda_1$ and $\lambda_2$ be two eigenvalues of $M$. If $Y_{\lambda_1},Y_{\lambda_2}$ are two associated eigenvectors   such that $[Y_{\lambda_1},Y_{\lambda_2}] \neq 0,$ then $\lambda_2=-1-\lambda_1$. In fact, 
$$[X_2,[Y_{\lambda_1},Y_{\lambda_2}] ]=(\lambda_1+\lambda_2)[Y_{\lambda_1},Y_{\lambda_2}] =a(\lambda_1+\lambda_2)X_1=-aX_1.$$
Since we assume $a\neq 0$, we obtain $-1=\lambda_1+\lambda_2.$ 
\item For any eigenvalue $\lambda$ of $M$ we consider an associated Jordan basis $\mathcal{B}_\lambda$ of the corresponding characteristic subvector space $C_\lambda$. The restriction to $M$ to  $C_\lambda$ is reduced in $\mathcal{B}_\lambda$ to a Jordan matrix composed of Jordan blocks $J_\lambda^s$, where $s$ denotes the dimensions of these blocks. If $\lambda \neq -\frac{1}{2}$, then $[C_\lambda,C_\lambda]=0.$ 
\item If $\lambda_1$ and $\lambda_2$ are two eigenvalues of $M$ with $\lambda_2 \neq -1-\lambda_1$, then $[C_{\lambda_1},C_{\lambda_2}]=0.$ 
\item 
If $\lambda \neq - \frac{1}{2}$ is an eigenvalue of $M$, 
then $-1-\lambda$ is also an eigenvalue of $M$ and the Segre characteristics (the ordered sequence of dimensions of Jordan blocks) of $\lambda$ and $-1-\lambda$ are equal. Moreover, if $J^s_{\lambda_1}$ and $J^r_{\lambda_2}$ are two Jordan blocks of $\lambda_1$ and $\lambda_2$ respectively with $s \neq r$, then $[J^s_{\lambda_1},J^r_{\lambda_2}]=0.$
\end{itemize}
All these remarks permit to reduce $M$ in the following form. Let $M_{2i+1,2i+2}^{2k+1,2k+2}$ the following submatrix of $M$ 
$$
\left(
\begin{array}{ll}
a^{2i+1}_{2k+1} & a^{2i+1}_{2k+2}\\
a^{2i+2}_{2k+1} & a^{2i+2}_{2k+2}\\
\end{array}
\right)
$$
Then we can reduce $M$ such that $M_{2i+1,2i+2}^{2k+1,2k+2}$ is equal to
$$
\left(
\begin{array}{ll}
0 & 1\\
0 & 0\\
\end{array}
\right)
$$
with $1 \leq i < k \leq p-1$ and $M_{2k+1,2k+2}^{2i+1,2i+2}=JM_{2i+1,2i+2}^{2k+1,2k+2}J$ with 
$$
J=\left(
\begin{array}{ll}
0 & 1\\
-1 & 0\\
\end{array}
\right)
$$
In particular the expression of $d\omega_1$ is preserved under this change of basis. \item In a second step, we contract any element of (\ref{7}) in an element of the family of the theorem but whose matrix $M$ is diagonal. We consider the Lie algebras (\ref{7}) whose Lie brackets are written in the Jordan basis defined in the first step. This means that if $\lambda$ is an eigenvalue of $M$ not equal to $\frac{1}{2}$, then $-1-\lambda$ is also eigenvalue and for any Jordan block of $\lambda$ of dimension $s$, there exists also a Jordan block of $-1-\lambda$ of dimension $s$ and if $\{Y_1,\cdots,Y_{s}\}$ and $\{Z_1,\cdots,Z_{s}\}$  are the Jordan basis of  these blocks, the Lie bracket of the Lie algebras (\ref{7}) are
\begin{equation}\label{8}
\left\{
\begin{array}{l}
\lbrack X_2,Y_1\rbrack=\lambda Y_1, \\
\lbrack X_2,Y_i\rbrack=\lambda Y_i+Y_{i-1},\ 2 \leq i \leq s,\\
\lbrack X_2,Z_s\rbrack=(-1-\lambda)Z_s, \\
\lbrack X_2,Z_i\rbrack=(-1-\lambda)Z_i-Y_{i+1},\ 1 \leq i \leq s-1.\\
\end{array}
\right.
\end{equation}
We consider the isomorphism $f_\varepsilon$ of $\C^{2p}$ given by
$$
\left\{
\begin{array}{l}
f_\varepsilon (X_1)=\varepsilon^{s-1}X_1,\\
f_\varepsilon (X_2)=X_2,\\
f_\varepsilon (Y_i)=\varepsilon^{i-1}Y_i,\ 1\leq i \leq s,\\
f_\varepsilon( Z_i)=\varepsilon^{s-i}Z_i, \ 1\leq i \leq s,\\
f_\varepsilon (U_l)=\varepsilon^{s-1}U_l,\\
f_\varepsilon (V_l)=V_l
\end{array}
\right.
$$
where $\{Y_i,Z_i,U_l,V_l\}$ is the Jordan basis of the step 1 with $d\omega_1(U_l,V_l)=1.$ This isomorphism defines a contraction whose associated matrix is diagonal in restriction to the space generated by the vectors $Y_i$ and $Z_i$. Moreover the expression of $d\omega_1$ is conserved. When we consider all the Jordan blocks, the contracted Lie algebra is of type (\ref{7}) but the matrix $M$ is diagonal. This proves the theorem. $\clubsuit$ 
\end{itemize}

\medskip

The parameters $\{a_1,\cdots,a_{p-1}\}$ are stable by contractions. The Lie algebra $\frak{f}$ corresponding to
$a_k=0$ satisfies
\[
\left\{
\begin{array}
[c]{l}%
d\omega_{1}=\omega_{1}\wedge\omega_{2}+\sum_{k=1}^{p-1}\omega_{2k+1}%
\wedge\omega_{2k+2},\\
d\omega_{2}=d\omega_{2k+1}=0,\;1\leq k\leq p-10,\\
d\omega_{2k+2}=- \omega_{2}\wedge\omega
_{2k+2},\;1\leq k\leq p-1.
\end{array}
\right.
\]
If $\psi_k$ is the $2$-cocycle of $\frak{f}$ such that $\psi_k(X_2,X_{2k+1})=X_{2k+1}$ and $\psi_k(X_2,X_{2k+2})=-X_{2k+2}$, then any Lie algebra of the Family \ref{frob} is a linear deformation of the Lie bracket $\mu_{\frak{f}}$ of $\frak{f}$ which writes $$\mu=\mu_{\frak{f}}+\sum a_k\psi_k.$$
In \cite{Giaquinto}, the notion of principal element of a Frobeniusian Lie algebra is defined. In the basis $\{X_1,\cdots,X_{2p}\}$ for which $d\omega_{1}=\omega_{1}\wedge\omega_{2}+\sum_{k=1}^{p-1}\omega_{2k+1}\wedge\omega_{2k+2}$, the principal element is $X_2$.
\begin{proposition}
The parameter $\{a_1,\cdots,a_{p-1}\}$ which are the invariants of Frobeniusian Lie algebras up a contraction are the eigenvalues of the principal element of $\g_{a_1,\cdots,a_{p-1}}$.
\end{proposition}

\medskip

\noindent{\bf Remarks.}
\begin{enumerate}
\item In \cite{GozeCras2}, we give the classification of frobeniusian Lie algebra up to a contraction when $\K=\R$.
\item Any algebras $\g_{a_1,\cdots,a_{p-1}}$ admit the following $\Z_2$-grading:
 $$\g_{a_1,\cdots,a_{p-1}}=\K\{X_1,X_{2}\} \oplus \K\{X_3,\cdots,X_{2p}\}.$$
 Thus the pairs $(\g_{a_1,\cdots,a_{p-1}},\K\{X_1,X_{2}\})$ is a symmetric pair with $\K\{X_1,X_{2}\}$ isomorphic to the affine $2$-dimensional Lie algebra.
\end{enumerate}
\bigskip

\section{Cartan class of $J$-invariant Pfaffian forms on a Lie group}
In this section, we are interested by the behavior of the Cartan class of non left invariant  Pfaffian forms on a Lie group $G$. We study two interesting cases, first the description of the $J$-invariant forms on the $3$-dimensional Heisenberg group, secondly a (non left invariant) contact form on the simple Lie group $SL(2n,\R)$. In this last case, we have seen that, if $n >1$, any left invariant Pfaffian form is not a contact form. Thus, it appears interesting to reduce the group of invariance to obtain a contact form.

\subsection{ $J$-invariant forms on the Heisenberg group }
The $3$-dimensional Heisenberg group $\HH_3$ is the real Lie group whose elements are the unipotent matrices
$$
\left(
\begin{array}{lll}
1 & x & z\\
0 & 1 & y \\
0 & 0 & 1
\end{array}
\right).
$$
We denote by $(x,y,z)$ the (global) coordinates of $\HH_3$ viewed as a $3$-dimensional differentiable manifold. Let $\T(\HH_3)$ the space of tangent vector fields on $\HH_3$, $\T_L(\HH_3)$ the subspace of left invariant vectors fields and $\T_R(\HH_3)$ the subspace  of right invariant vectors fields. Similarly, we denote by $\bigwedge^p(\HH_3)$ (respectively $\bigwedge^1_L(\HH_3)$) the space of left $p$-forms on  $\HH_3$ (respectively left invariant $p$-forms).
In any point $g=(x,y,z)$ of $\HH_3$, the linear space $\T_L(\HH_3)$ is of dimension $3$ and generated by the vector fields
$$
\displaystyle
\left\{
\begin{array}{l}
  X_1= \frac{\partial}{\partial x},\\
  X_2= \frac{\partial}{\partial y}+x \frac{\partial}{\partial z},\\
  X_3=\frac{\partial}{\partial z}.
\end{array}
\right.
$$
We have $[X_1,X_2]=X_3, \ [X_1,X_3]=[X_2,X_3]=0.$ These relations determine the structure of Lie algebra of the $3$-dimensional Heisenberg algebra $\h_3$, whose corresponding elements will be also denoted by $X_1,X_2,X_3$. In any point $g=(x,y,z)$ of $\HH_3$, the linear space $\T_R(\HH_3)$ is of dimension $3$ and generated by the vector fields
$$
\displaystyle
\left\{
\begin{array}{l}
  \widetilde{X_1}= \frac{\partial}{\partial x}+y \frac{\partial}{\partial z},\\
  \widetilde{X_2}= \frac{\partial}{\partial y},\\
  \widetilde{X_3}=\frac{\partial}{\partial z}.
\end{array}
\right.
$$
These vector fields generate a $3$-dimensional Lie algebra, $\widetilde{\h_3}$, isomorphic to $\h_3$ and whose Lie bracket is given by $[ \widetilde{X_1}, \widetilde{X_2}]=- \widetilde{X_3}.$ A basis of $\bigwedge^1_L(\HH_3)$, identified with the dual space $\h_3^*$, is given by
$$
\displaystyle
\left\{
\begin{array}{l}
  \omega_1= dx,\\
 \omega_2= dy,\\
  \omega_3=dz-xdy.
\end{array}
\right.
$$
We easily verify  the Maurer-Cartan equations, $d\omega_1=d\omega_2=0,d\omega_3=-\omega_1\wedge \omega_2.$ The left invariance of the Pfaffian forms $\omega_i$ is equivalent to $L_{\widetilde{U}}\omega_i=0$ for any $U \in \widetilde{\h_3}$.

Every $1$-dimensional Lie subalgebra of $\widetilde{\h_3}$ is generated by a non null vector of $\widetilde{\h_3}$. For $2$-dimensional Lie  subalgebras, we have
\begin{lemma}
Any $2$-dimensional Lie subalgebra of $\widetilde{\h_3}$ is generated by the vectors
$$U_1=\alpha_1\widetilde{X_1}+\alpha_2\widetilde{X_2}, \ U_2=\widetilde{X_3}$$
with $\alpha_1^2+\alpha_2^2 \neq 0$.
\end{lemma}
\pf In fact, since $\widetilde{\h_3}$ is nilpotent, any Lie subalgebra is nilpotent and, in this case, abelian. If $\{U_1,U_2\}$ is a basis, thus $[U_1,U_2]=0$. If $U_i=\sum a_{ij}\widetilde{X_j}$, thus $a_{11}a_{22}-a_{12}a_{21}=0.$
We can take $U_2=\widetilde{X_3}$ and $U_1=a_{11}\widetilde{X_1}+a_{12}\widetilde{X_2}.$

Let $J$ the subalgebra of $\widetilde{\h_3}$ generated by $U_1$ and $U_2$.
\begin{proposition}
Any $J$-invariant Pfaffian form $\omega \in \Lambda^1(\HH_3)$ is written $w=a_1\omega_1+a_2\omega_2+a_3\omega_3$ with
$$
\displaystyle
\left\{
\begin{array}{ll}
\medskip
 \displaystyle \alpha_1\frac{\partial a_1}{\partial x}+\alpha_2\frac{\partial a_1}{\partial y}=0, & \displaystyle\frac{\partial a_1}{\partial z}=0,\\
  \medskip
\displaystyle \alpha_1\frac{\partial a_2}{\partial x}+\alpha_2\frac{\partial a_2}{\partial y}=0, &\displaystyle \frac{\partial a_2}{\partial z}=0,\\
  \medskip
 \displaystyle \alpha_1\frac{\partial a_3}{\partial x}+\alpha_2\frac{\partial a_3}{\partial y}=0, & \displaystyle\frac{\partial a_3}{\partial z}=0.\\
\end{array}
\right.
$$
\end{proposition}
\pf In fact the equation $L_{U_2}\omega=L_{\widetilde{X_3}}=0$ implies that $\displaystyle \frac{\partial a_i}{\partial z}=0$ for $i=1,2,3$. Thus the functions $a_i(x,y,z)$ only depend  on $x$ and $y$. The second equation $L_{U_2}\omega=0$ implies $\displaystyle \alpha_1\frac{\partial a_i}{\partial x}+\alpha_2\frac{\partial a_i}{\partial y}=0.$

\medskip

\noindent{\bf Remark.} Since $\alpha_1$ and $\alpha_2$ are not simultaneously equal to $0$, the equations $\displaystyle \alpha_1\frac{\partial a_i}{\partial x}+\alpha_2\frac{\partial a_i}{\partial y}=0$ correspond to the transport equations. Assume that $\alpha_1 \neq 0$. If we put $\alpha=\frac{\alpha_2}{\alpha_1}$ and $b_i(y)=a_i(0,y)$, we have as unique solution
$$a_i(x,y)=b_i(y-\alpha x).$$
The form $\omega$ writes $\omega=a_1dx+(a_2-a_3x)dy+a_3dz.$ Putting $y'=y-\alpha x$, thus $dy'=dy - \alpha dx$ and
$$\omega= (b_1(y')+\alpha b_2(y') -\alpha xb_3(y'))dx+(b_2(y')-xb_3(y'))dy' +b_3(y')dz.$$
This form is then of type
$$\omega = (c_1(v)-\alpha uc_2(v))du+(c_2(v)-uc_3(v))dv+c_3(v)dw.$$

\medskip

Let $\Sigma=\{g \in \HH_3, \ \omega_g(U_g)=0 \ \forall U \in J\}$ be the singular set associated with $\omega$. To simplify notations, we assume that $alpha_1 \neq 0$ and $ \alpha=\alpha_2\alpha_1^{-1}$. In this case
$U=a(\widetilde{X_1}+\alpha \widetilde{X_2})+b\widetilde{X_3}.$ Since $\omega=b_1(y-\alpha x)\omega_1 +b_2(y-\alpha x)\omega_2+b_3(y-\alpha x)\omega_3,$ thus $\omega_g(\widetilde{X_3})_g=0$ implies $b_3(y-\alpha x)=0$ with $g=(x,y,z)$. Likewise $\omega(\widetilde{X_1}+\alpha \widetilde{X_2})=b_1(y-\alpha x) +\alpha b_2(y-\alpha x)-\alpha xb_3(y-\alpha x)=0.$
\begin{proposition}
The singular set $\Sigma$ associated with the $J$-invariant form $\omega$ is the set of $g=(x,y,z) \in \HH_3$ such that
$$
\left\{
\begin{array}{l}
b_3(y-\alpha x)=0,\\
b_1(y-\alpha x) +\alpha b_2(y-\alpha x)=0.
\end{array}
\right.
$$
\end{proposition}

\medskip

Assume now that the $J$-invariant form $\omega$ is a {\bf contact form} on $\HH_3$, that is $\omega \wedge d\omega \neq 0$ in any point of $\HH_3$. Since the functions $b_i$ can be considered as functions of one variable, we will denote by $b_i'$ its derivative.  The condition $\omega \wedge d\omega \neq 0$ is then equivalent to
$$b_1b'_3-b'_1b_3+\alpha (b_2b'_3-b'_2b_3)-b_3^2 \neq 0.$$
\begin{proposition}
If $\omega$ is a $J$-invariant contact form on $\HH_3$, then its singular set $\Sigma$  is empty.
\end{proposition}
\pf If $\g=(x,y,z) \in \Sigma$, then $b_3(y-\alpha x)=0$. In this case
$$\omega \wedge d\omega (g)= (b_1(y-\alpha x) +\alpha b_2(y-\alpha x))b'_3(y-\alpha x)dx\wedge dy\wedge dz.$$
But, if $ g \in \Sigma$, $b_1(y-\alpha x) +\alpha b_2(y-\alpha x)=0$ and $\omega \wedge d\omega (g)=0$. This contradicts the hypothesis.

\subsection{A contact form on the simple group $SL(2n,\R)$}

Recall that $SL(2n,\R)$ is the simple Lie group constitued of $2n \times 2n$ real matrices with determinant  equal to $1.$
Let $M=(x_{ij})$ be in $SL(2n,\R).$ We denote by $X_{ij}$ the minor of the coefficient $x_{ij}$ and by $det M$ the determinant of $M.$

\begin{lemma}
We have 
$$\Delta= d(det M)=\sum_{1\leq i \leq j \leq n} (-1)^{i+j}X_{ij}dx_{ij}$$
\end{lemma}
 
We consider the Pfaffian form $\tilde{\omega}$ on $GL(2n,\R)$
$$\tilde{\omega}=\sum_{ j=1}^{2n} \sum_{  i=1}^{ n} x_{j,2i-1}dx_{j,2i}-x_{j,2i}dx_{j,2i-1}$$
and $\omega =i^{*} \tilde{\omega}$ its induced form on $SL(2n,\R)$ when $i: SL(2n,\R)\rightarrow GL(2n,\R)$ is the natural injection.

\begin{proposition}
The Pfaffian form $\omega$ is a contact form on $SL(2n,\R).$
\end{proposition}

\pf The form $\omega$ is a contact form on $SL(2n,\R)$ if and only if 
$$\tilde{\omega} \wedge (d\tilde{\omega})^{2n-1}\wedge d\Delta=-2^{2n^2}n\Delta.$$
Since $\Delta=1$ on $SL(2n,\R),$ $\omega$ is a contact form.

\begin{corollary}
The Lie group $SL(2n,\R)$ is a contact manifold.
\end{corollary}

The Lie group $SO(2n,\R)=\{M \in GL(2n,\R) / ^t \! M\cdot M=Id\}$ is a compact Lie subgroup of $SL(2n,\R).$  
\begin{theorem}
The Pfaffian form $\omega$ is a contact form on $SL(2n,\R)$ invariant by $SO(2n,\R).$
\end{theorem}

\pf The Reeb vector field $Z_\omega$ associated with the contact form $\omega$ is  the vector field induced on $SL(2n,\R)$ by $$\widetilde{\displaystyle Z_{\omega}}=\frac{1}{2n}\sum_{j=1}^{2n} \sum_{i=1}^n (-1)^{j+1} \left(X_{j,2i-1}\frac{\partial}{\partial x_{j,2i}}+X_{j,2i}\frac{\partial}{\partial x_{j,2i-1}}\right).$$
In fact $\omega (Z_{\omega})=\frac{1}{2n}(2n\Delta)=1$ and $i(Z_{\omega})(d\omega)=d\Delta=0.$
Let us show that $Z_\omega$ is invariant by $SO(2n,\R).$ For any $M\in GL(2n,\R),$ the jacobian matrix at the unity $e=I_{2n}$ of the left translation by $M$ is $M\otimes I_{2n}.$ Then 
$(L_M)_{e}^T(Z_\omega)(e)=(Z_\omega)(M)$ if and only if the matrix $M=(x_{ij})$ satisfies 
$X_{ij}=(-1)^{i+j}x_{ij}.$ We deduce that $M \cdot \, ^t \! M$ and $M \in SO(2n,\R).$

\begin{proposition}
The singular set 
$\Sigma_\omega$ of the $SO(2n,\R)$-invariant contact form $\omega$ on $SL(2n,\R)$ is empty.
\end{proposition}

\pf A basis of the Lie algebra of $SL(2n,\R)$ is given by the field 
$i^{*}(A_{ij})$ where $i$ is the injection of $SL(2n,\R)$ in $GL(2n,\R)$ and 
$$\displaystyle A_{ij}=\sum_{l=1}^{2n}  \left( x_{jl}\frac{\partial}{\partial x_{il}}-x_{il}\frac{\partial}{\partial x_{jl}}\right).$$

Then $\tilde{\omega}(A_{ij})=\sum_{l=1}^{2n} (x_{i,2l-1}x_{j,2l}-x_{i,2l}x_{j,2l-1}) .$

The singular set is the set of the points $M\in SL(2n,\R)$ such that $\tilde{\omega}(A_{ij})(M)=0$ for all $i,j \in \{1,\cdots , 2n \}$.
We deduce that $\Sigma_\omega$ is the real algebraic subvariety of $SL(2n,\R)$ given by the system of polynomial equations given by 
$$ \sum_{l=1}^{2n} (x_{i,2l-1}x_{j,2l}-x_{i,2l}x_{j,2l-1})$$ 
for $i,j \in \{1,\cdots , 2n \}$.

\end{document}